\documentclass[12pt,twoside]{article}
\usepackage[english]{babel}
\usepackage[latin1]{inputenc}
\usepackage{amsmath}
\usepackage{amssymb,amsfonts}
\usepackage{graphicx}                   

\newcommand{\bE}{\mathbf{E}}
\newcommand{\bG}{\mathbf{G}}
\newcommand{\bH}{\mathbf{H}}

\newcommand{\bS}{\mathbf{S}}

\newcommand{\ba}{\mathbf{a}}
\newcommand{\bc}{\mathbf{c}}
\newcommand{\be}{\mathbf{e}}

\newcommand{\bl}{\mathbf{l}}
\newcommand{\bh}{\mathbf{h}}
\newcommand{\bj}{\mathbf{j}}
\newcommand{\bq}{\mathbf{q}}
\newcommand{\bk}{\mathbf{k}}

\newcommand{\EUC}{\bE^3}

\newcommand{\HYN}{\bH^n}
\newcommand{\HYP}{\bH^3}

\begin{document}
\pagestyle{myheadings}
\markboth{\centerline{Emil Moln\'ar and Jen\H o Szirmai}}
{Top dense hyperbolic ball packings and coverings \dots}
\title
{Top dense hyperbolic ball packings and coverings for complete Coxeter orthoscheme groups}

\author{\normalsize{Emil Moln\'ar and Jen\H o Szirmai} \thanks{} \\
\normalsize Budapest University of Technology and \\
\normalsize Economics Institute of Mathematics, \\
\normalsize Department of Geometry \\
\date{\normalsize{\today}}}


\maketitle


\begin{abstract}
In $n$-dimensional hyperbolic space $\mathbf{H}^n$ $(n\ge2)$ there are $3$-types of spheres (balls):
the sphere, horosphere and hypersphere. If $n=2,3$ we know an universal upper bound of the ball packing densities,
where each ball volume is related to the volume of the corresponding Dirichlet-Voronoi (D-V) cell. E.g.
in $\HYP$ the densest horoball packing is derived from the
$\{3,3,6\}$ Coxeter tiling consisting of ideal regular simplices $T_{reg}^\infty$  with dihedral angles
$\frac{\pi}{3}$. The density of this packing is $\delta_3^\infty\approx 0.85328$
and this provides a very rough upper bound for the ball packing densities as well. However,
there are no "essential" results regarding the "classical" ball packings with congruent balls, and for ball coverings either.

The goal of this paper to find
the extremal ball arrangements in $\HYP$ with "classical balls". We consider only periodic congruent ball arrangements (for simplicity)
related to the generalized, so-called {\it complete Coxeter orthoschemes} and their extended groups.
In Theorems 1.1-1.2 we formulate also conjectures for the densest ball packing with density $0.77147\dots$ and the loosest
ball covering with density $1.36893\dots$, respectively.
Both are related with the extended Coxeter group $(5, 3, 5)$ and the so-called hyperbolic football manifold (look at Fig.~3).
These facts can have important relations with fullerens in crystallography.
\end{abstract}

 \newtheorem{thm}{Theorem}[section]
 \newtheorem{prop}{Proposition}[section]
 \newtheorem{lem}{Lemma}[section]
 \newtheorem{cor}{Corollary}[section]
 \newtheorem{exm}{Example}[section]
 \newtheorem{dfn}{Definition}[section]
 \newtheorem{rem}{Remark}[section]
 \numberwithin{equation}{section}

\maketitle

\section{Introduction}
Ball packing problems concern the arrangements of non-overlapping equal balls which fill space.
Usually, space is the classical three-dimensional Euclidean space $\EUC$. However, ball packing problems
can be generalized to the other $3$-dimensional Thurston geometries.

In an $n$-dimensional space of constant curvature $\bE^n$, $\bH^n$, $\bS^n$ $(n\ge2)$ let $d_n(r)$
be the density of $n+1$ equal balls of radius $r$ mutually touching one another with respect to the simplex
spanned by the centres of the balls. {L.~Fejes T\'oth} and {H.~S.~M.~Coxeter} conjectured that in an
$n$-dimensional space of constant curvature the density of packing balls of radius $r$ cannot exceed $d_n(r)$.
This conjecture has been proved by C.~Rogers in Euclidean $n$-space.
The 2-dimensional case has been solved by {L.~Fejes T\'oth}.
In a three-dimensional space of constant curvature the problem has been investigated
by { K.~B\"or\"oczky} and {A.~Florian} in \cite{BF64} and it has been studied by
{ K.~B\"or\"oczky} in \cite{B78} for $n$-dimensional spaces of constant curvature $(n\ge 4)$.
The upper bound $d_n(\infty)$ for $\HYN$ $(n=2,3)$ is attained for a regular horoball packing, that is, a
packing by horoballs which are inscribed in the cells of a regular honeycomb of $\overline{\mathbf{H}}^n$
(i.e. $\HYN$ is closed by its ideal points, or ends).

In $\overline{\mathbf{H}}^3$ there is exactly one horoball packing with horoballs in the same type whose
Dirichlet--Voronoi cells give rise to a
regular honeycomb described by the Schl\"afli symbol $\{6,~3,~3\}$ . Its
dual $\{3,3,6\}$ consists of ideal regular simplices $T_{reg}^\infty$  with dihedral angles $\frac{\pi}{3}$ building
up a 6-cycle around each edge
of the tessellation. The density of this packing is $\delta_3^\infty\approx 0.85328$.
We have considered some new aspects related to the horoball and hyperball packings in
\cite{Sz12-1}, \cite{Sz13-1}, \cite{Sz07-1}, \cite{KSz2} and \cite{KSz1}.

However, there are no "essential" results regarding the "classical ball packings and coverings" with congruent balls.
What are the extremal ball arrangements in $\mathbf{H}^n$ and what are their densities?

The goal of this paper to study the above problems in $\HYP$ with "classical balls". We consider
periodic congruent ball packings and coverings (for simplicity and for good constructions, only)
related to the generalized, so-called complete Coxeter orthoschemes and their extended groups. We formulate two
theorems and conjectures for the densest ball packing with density $0.77147\dots$, and for the loosest ball covering with
density $1.36893\dots$, respectively.
\begin{dfn} For a given packing of $\HYN$ ($2\le n \in \mathbf{N}$ natural numbers) each ball volume is related to the volume of its D-V cell,
(i.e. the closed domain whose any point lies not further to the given ball centre than to each other centre of the ball system),
then take the infimum of these ratios for all balls; this is called the density of the given
packing. Then the supremum of these infima is taken for all ball packings of $\HYN$ to get the maximal packing density and the
densest packing of $\HYN$ (if such a packing exists at all).
\end{dfn}
\begin{thm}
Among the extended complete Coxeter orthoscheme groups the group $(5, 3, 5)$, extended with a half-turn, provides the ball
centre orbit of $A_3$ for the densest ball packing by its football D-V-cell $\{5, 6, 6\}$ and inscribed ball,
see Fig.~1 and 3. The density of this ball packing is $0.77147 \dots$ with radius $r = 0.95142... = A_3F_{03}$ .
This is the conjectured maximal density for all ball packings of $\HYP$.
\end{thm}
\begin{dfn}
For a given covering of $\HYN$ each ball volume is related to the volume of its D-V cell, than take the supremum of
these ratios for all balls (this is called the density of the given covering); then the infimum
of these suprema is taken for all ball coverings of $\HYN$ to get the minimal covering density and the loosest covering of $\HYN$
(if such a covering exists at all).
\end{dfn}
\begin{thm}
Among the above groups the same (with half-turn extended $(5, 3, 5)$ group) provides the same
(congruent) ball centre orbit of $A_3$ for the loosest ball covering, again by its football
D-V-cell $\{5, 6, 6\}$ but with its circumscribed ball, see Fig.~1 and 3.
The density of this ball covering is $1.36893 \dots$ with radius $R = 1.12484... = A_3F_{12}$.
This is the conjectured minimal density for all ball coverings of $\HYP$.
\end{thm}
Our numerical results are collected in Tables, where we give only the relevant data (because of the page restriction of this paper).

{\bf Our systematic computations will give the Proof of Theorems 1.1-1.2.}

We will use the well-known Beltrami-Cayley-Klein model of $\HYP$ with the classical
projective metric calculus (see e.g. \cite{M97}).
\section{The projective model and complete orthoschemes}
We use for $\mathbf{H}^3$ (and analogously for $\HYN$, $n>3$) the projective model in the Lorentz space $\mathbf{E}^{1,3}$
that denotes the real vector space $\mathbf{V}^{4}$ equipped with the bilinear
form of signature $(1,3)$,
$
\langle \mathbf{x},~\mathbf{y} \rangle = -x^0y^0+x^1y^1+x^2y^2+ x^3 y^3,
$
where the non-zero vectors
$
\mathbf{x}=(x^0,x^1,x^2,x^3)\in\mathbf{V}^{4} \ \  \text{and} \ \ \mathbf{y}=(y^0,y^1,y^2,y^3)\in\mathbf{V}^{4},
$
are determined up to real factors, for representing points of $\mathcal{P}^n(\mathbf{R})$. Then $\mathbf{H}^3$ can be interpreted
as the interior of the quadric
$
Q=\{(\mathbf{x})\in\mathcal{P}^3 | \langle  \mathbf{x},~\mathbf{x} \rangle =0 \}=:\partial \mathbf{H}^3
$
in the real projective space $\mathcal{P}^3(\mathbf{V}^{4},
\mbox{\boldmath$V$}\!_{4})$ (here $\mbox{\boldmath$V$}\!_{4}$ is the dual space of $\mathbf{V}^{4}$).
Namely, for an interior point $Y(\mathbf{y})$ holds $\langle  \mathbf{y},~\mathbf{y} \rangle <0$.

The points of the boundary $\partial \mathbf{H}^3 $ in $\mathcal{P}^3$
are called points at infinity, or at the absolute of $\mathbf{H}^3 $. The points lying outside $\partial \mathbf{H}^3 $
are said to be outer points of $\mathbf{H}^3 $ relative to $Q$. Let $(\mathbf{x}) \in \mathcal{P}^3$, a point
$(\mathbf{y}) \in \mathcal{P}^3$ is said to be conjugate to $(\mathbf{x})$ relative to $Q$ if
$\langle \mathbf{x},~\mathbf{y} \rangle =0$ holds. The set of all points which are conjugate to $(\mathbf{x})$
form a (polar) hyperplane
$
pol(\mathbf{x}):=\{(\mathbf{y})\in\mathcal{P}^3 | \langle \mathbf{x},~\mathbf{y} \rangle =0 \}.
$
Thus the quadric $Q$ induces a bijection
(linear polarity $\mathbf{V}^{4} \rightarrow
\mbox{\boldmath$V$}\!_{4})$
from the points of $\mathcal{P}^3$ onto their polar hyperplanes.

The point $X (\bold{x})$ and the hyperplane $\alpha (\mbox{\boldmath$a$})$
are incident if $\bold{x}\mbox{\boldmath$a$}=0$ ($\bold{x} \in \bold{V}^{4} \setminus \{{\mathbf{0}}\}, \ \mbox{\boldmath$a$}
\in \mbox{\boldmath$V$}_{4}
\setminus \{\mbox{\boldmath{$0$}}\}$).
The constant $k =\sqrt{\frac{-1}{K}}$ is the natural length unit in
$\mathbf{H}^3$, where $K$ denotes the constant negative sectional curvature. In the following we may assume that $k=1$.
\subsection{Characteristic orthoschemes and their volumes}

{\it An orthoscheme $\mathcal{O}$ in $\mathbf{H}^n$ $n \geq 2$ in classical sense}
is a simplex bounded by $n+1$ hyperplanes $H_0,\dots,H_n$
such that $H_i \bot H_j, \  \text{for} \ j \ne i-1,i,i+1.$ Or,
equivalently, the $n+1$ vertices of $\mathcal{O}$ can be
labelled by $A_0,A_1,\dots,A_n$ in such a way that
$\text{span}(A_0,\dots,A_i) \perp \text{span}(A_i,\dots,A_n) \ \ \text{for} \ \ 0<i<n-1.$

Geometrically, complete orthoschemes of degree $m=0,1,2$ can be described as
follows:

\begin{enumerate}
\item
For $m=0$, they coincide with the class of classical orthoschemes introduced by
{{Schl\"afli}}. The initial and final vertices, $A_0$ and $A_n$ of the orthogonal edge-path
$A_iA_{i+1},~ i=0,\dots,n-1$, are called principal vertices of the orthoscheme (see Remark 4.1).
\item
A complete orthoscheme of degree $m=1$ can be constructed from an
orthoscheme with one outer principal vertex, one of $A_0$ or $A_n$, which is simply truncated by
its polar plane (see Fig.~1-2).
\item
A complete orthoscheme of degree $m=2$ can be constructed from an
orthoscheme with two outer principal vertices, $A_0$ and $A_n$, which is doubly truncated by
their polar planes $pol(A_0)$ and $pol(A_n)$ (see Fig.~1-2).
\end{enumerate}
For the {\it complete Coxeter orthoschemes} $\mathcal{O} \subset \mathbf{H}^n$ we adopt the usual
conventions and sometimes even use them in the Coxeter case: If two nodes are related by the weight $\cos{\frac{\pi}{p}}$
then they are joined by a ($p-2$)-fold line for $p=3,~4$ and by a single line marked by $p$ for $p \geq 5$.
In the hyperbolic case if two bounding hyperplanes of $S$ are parallel, then the corresponding nodes
are joined by a line marked $\infty$. If they are divergent then their nodes are joined by a dotted line.

In the following we concentrate only on dimensions $3$ and on hyperbolic
Coxeter-Schl$\ddot{a}$fli symbol of the complete orthoscheme tiling $\mathcal{P}$ generated by a complete orthoscheme $\mathcal{O}$.
To every scheme there is a corresponding
symmetric $4 \times 4$ matrix $(b^{ij})$ where $b^{ii}=1$ and, for $i \ne j\in \{0,1,2,3\}$,
$b^{ij}$ equals $-\cos{\alpha_{ij}}$ with all angles $\alpha_{ij}$ between the faces $i$,$j$ of $\mathcal{O}$.

For example, $(b^{ij})$ below is the so called Coxeter-Schl\"afli matrix with
parameters $(u;v;w)$, i.e. $\alpha_{01}=\frac{\pi}{u}$, $\alpha_{12}=\frac{\pi}{v}$, $\alpha_{23}=\frac{\pi}{w}$
to be discussed yet for hyperbolicity. Now only $3\le u,v,w$ come into account (see \cite{IH85}, \cite{IH90}). Then it holds
\[
(b^{ij})=\langle \mbox{\boldmath$b^i$},\mbox{\boldmath$b^j$} \rangle:=\begin{pmatrix}
1& -\cos{\frac{\pi}{u}}& 0 & 0 \\
-\cos{\frac{\pi}{u}} & 1 & -\cos{\frac{\pi}{v}}& 0 \\
0 & -\cos{\frac{\pi}{v}} & 1 & -\cos{\frac{\pi}{w}} \\
0 & 0 & -\cos{\frac{\pi}{w}} & 1
\end{pmatrix}. \tag{2.1}
\]

This $3$-dimensional complete (truncated or frustum) orthoschemes $\mathcal{O}=W_{uvw}$ and its reflection group $\bG_{uvw}$ wil be described in
Fig.~2, and by the symmetric Coxeter-Schl\"afli matrix $(b^{ij})$ in formula (2.1), furthermore by its inverse matrix $(a_{ij})$ in formula
(2.2).
\[
\begin{gathered}
(a_{ij})=(b^{ij})^{-1}=\langle \ba_i, \ba_j \rangle:=\\
=\frac{1}{B} \begin{pmatrix}
\sin^2{\frac{\pi}{w}}-\cos^2{\frac{\pi}{v}}& \cos{\frac{\pi}{u}}\sin^2{\frac{\pi}{w}}& \cos{\frac{\pi}{u}}\cos{\frac{\pi}{v}} & \cos{\frac{\pi}{u}}\cos{\frac{\pi}{v}}\cos{\frac{\pi}{w}} \\
\cos{\frac{\pi}{u}}\sin^2{\frac{\pi}{w}} & \sin^2{\frac{\pi}{w}} & \cos{\frac{\pi}{v}}& \cos{\frac{\pi}{w}}\cos{\frac{\pi}{v}} \\
\cos{\frac{\pi}{u}}\cos{\frac{\pi}{v}} & \cos{\frac{\pi}{v}} & \sin^2{\frac{\pi}{u}}  & \cos{\frac{\pi}{w}}\sin^2{\frac{\pi}{u}}  \\
\cos{\frac{\pi}{u}}\cos{\frac{\pi}{v}}\cos{\frac{\pi}{w}}  & \cos{\frac{\pi}{w}}\cos{\frac{\pi}{v}} & \cos{\frac{\pi}{w}}\sin^2{\frac{\pi}{u}}  & \sin^2{\frac{\pi}{u}}-\cos^2{\frac{\pi}{v}}
\end{pmatrix}, \tag{2.2}
\end{gathered}
\]
where
$$
B=\det(b^{ij})=\sin^2{\frac{\pi}{u}}\sin^2{\frac{\pi}{w}}-\cos^2{\frac{\pi}{v}} <0, \ \ \text{i.e.} \ \sin{\frac{\pi}{u}}\sin{\frac{\pi}{w}}-\cos{\frac{\pi}{v}}<0.
$$

In the following we use the above orthoschemes whose volume is derived by the next
Theorem of {{R.~Kellerhals}} (\cite{K89}, by the ideas of N.~I.~Lobachevsky):
\begin{thm}{\rm{(R.~Kellerhals)}} The volume of a three-dimensional hyperbolic
complete orthoscheme $\mathcal{O}=W_{uvw} \subset \mathbf{H}^3$
is expressed with the essential
angles $\alpha_{01}=\frac{\pi}{u}$, $\alpha_{12}=\frac{\pi}{v}$, $\alpha_{23}=\frac{\pi}{w}$, $(0 \le \alpha_{ij}
\le \frac{\pi}{2})$
(Fig.~1.a,~b) in the following form:

\begin{align}
&\mathrm{Vol}(\mathcal{O})=\frac{1}{4} \{ \mathcal{L}(\alpha_{01}+\theta)-
\mathcal{L}(\alpha_{01}-\theta)+\mathcal{L}(\frac{\pi}{2}+\alpha_{12}-\theta)+ \notag \\
&+\mathcal{L}(\frac{\pi}{2}-\alpha_{12}-\theta)+\mathcal{L}(\alpha_{23}+\theta)-
\mathcal{L}(\alpha_{23}-\theta)+2\mathcal{L}(\frac{\pi}{2}-\theta) \}, \notag
\end{align}
where $\theta \in [0,\frac{\pi}{2})$ is defined by:
$$
\tan(\theta)=\frac{\sqrt{ \cos^2{\alpha_{12}}-\sin^2{\alpha_{01}} \sin^2{\alpha_{23}
}}} {\cos{\alpha_{01}}\cos{\alpha_{23}}},
$$
and where $\mathcal{L}(x):=-\int\limits_0^x \log \vert {2\sin{t}} \vert dt$ \ denotes the
Lobachevsky function.
\end{thm}
The volume $\mathrm{Vol}(B(R))$ of a {\it ball} $B(R)$ of radius $R$ can be computed by the classical formula of J.~Bolyai:
\begin{equation}
\begin{gathered}
\mathrm{Vol}(B(R))=2\pi (\cosh(R) \sinh(R) -R)=\pi(\sinh(2R)-2R)=\\
=\frac{4}{3} \pi R^3(1+\frac{1}{5}R^2+\frac{2}{105}R^4+\dots). \tag{2.3}
\end{gathered}
\end{equation}

\section{Essential points in a complete (truncated) orthoscheme}
Let $A_0(\ba_0)$, $A_1(\ba_1)$, $A_2(\ba_2)$, $A_3(\ba_3)$ be the vertices of the above complete orthoscheme $W_{uvw}$ by the
$u,v,w$ above (see Fig.~1,2).
The principal vertices $A_0$ and $A_3$ can be proper ($a_{ii}<0$, $i \in \{0,3\})$, boundary ($a_{ii}=0$) or outer points $(a_{ii}>0)$ We exploit the logical symmetry $0,1 \leftrightarrow 3,2$.

We distinguish the following main configurations of the principal verices $A_0$ and $A_3$:
\begin{enumerate}
\item[{\bf 1.}] {\it $A_3$ is proper or boundary point $\frac{\pi}{u}+\frac{\pi}{v} \ge \frac{\pi}{2}$.}
\begin{enumerate}
\item[{\bf 1.i.}] $A_0$ is proper or boundary point $\frac{\pi}{v}+\frac{\pi}{w} \ge \frac{\pi}{2}$.
\item[{\bf 1.s.i}] $u=w$, $F_{03}F_{12}$ is half turn axis, $\mathbf{h}$ is the half turn changing $0 \leftrightarrow 3$, $1 \leftrightarrow 2$.
Here a "half orthoscheme" $JQEB_{13}F_{12}B_{02}F_{03}A_2$ will be the fundamental domain of $\bG_{u=w,v}$.
\item[{\bf 1.ii}] $A_0$ is outer $\frac{\pi}{v}+\frac{\pi}{w} < \frac{\pi}{2}$, then $a_0(\mbox{\boldmath$a$}_{0})=CLH$ is its polar plane.
\end{enumerate}
\item[{\bf 2.}] {\it $A_3$ is outer point $\frac{\pi}{u}+\frac{\pi}{v} < \frac{\pi}{2}$, then $a_3(\mbox{\boldmath$a$}_{3})=JEQ$ is its polar plane.}
\begin{enumerate}
\item[{\bf 2.i.}] $A_0$ is proper or boundary point $\frac{\pi}{v}+\frac{\pi}{w} \ge \frac{\pi}{2}$.
\item[{\bf 2.ii}] $A_0$ is also outer $\frac{\pi}{v}+\frac{\pi}{w} < \frac{\pi}{2}$, then $a_0(\mbox{\boldmath$a$}_{0})=CLH$ is its polar plane.
\item[{\bf 2.s.ii}] $u=w$, $F_{03}F_{12}$ is half turn axis, $\mathbf{h}$ is the half turn changing $0 \leftrightarrow 3$, $1 \leftrightarrow 2$.
Here a "half orthoscheme" $JQEB_{13}F_{12}B_{02}F_{03}A_2$ will be the fundamental domain of $\bG_{u=w,v}$.
\end{enumerate}
\end{enumerate}
We obtain with easy calculations the following inportant lemmas (see Fig.~1,2):
\begin{lem}
Let $A_0$ be an outer principal vertex of the orthoscheme $W_{uvw}$ and let $a_0(\mbox{\boldmath$a$}_{0})=CLH$ be its polar plane where
$C=a_0 \cap A_0A_1$, $L=a_0 \cap A_0A_2$, $H=a_0 \cap A_0A_3$ whose vectors are the following:
\begin{equation}
\begin{gathered}
C(\bc)=a_0 \cap A_0A_1; ~ \bc=\ba_1-\frac{a_{01}}{a_{00}} \ba_0, \ \mathrm{with} \\
\langle \bc,\bc \rangle =\frac{(a_{11}a_{00}-a_{01}^2)}{a_{00}}=\langle \bc,\ba_1 \rangle
=\frac{\sin^2\frac{\pi}{w}}{\sin^2\frac{\pi}{w}-\cos^2\frac{\pi}{v}}=\frac{a_{11}}{a_{00}} \\
L(\mathbf{l})=a_0 \cap A_0A_2; ~ \mathbf{l}=\ba_2-\frac{a_{02}}{a_{00}} \ba_0, \ \mathrm{with} \\
\langle \bl,\bl \rangle =\frac{(a_{22}a_{00}-a_{02}^2)}{a_{00}}=\langle \bl,\ba_2 \rangle
=\frac{1}{\sin^2\frac{\pi}{w}-\cos^2\frac{\pi}{v}}=\frac{1}{Ba_{00}} \\
H(\bh)=a_0 \cap A_0A_3; ~ \bh=\ba_3-\frac{a_{03}}{a_{00}} \ba_0, \ \mathrm{with} \\
\langle \bh,\bh \rangle =\frac{(a_{33}a_{00}-a_{03}^2)}{a_{00}}=\langle \bh,\ba_3 \rangle
=\frac{\sin^2\frac{\pi}{v}}{\sin^2\frac{\pi}{w}-\cos^2\frac{\pi}{v}}=\frac{\sin^2\frac{\pi}{v}}{Ba_{00}}
\end{gathered} \tag{3.1}
\end{equation}
\end{lem}
\begin{figure}[htb]
\centering
\includegraphics[width=60mm]{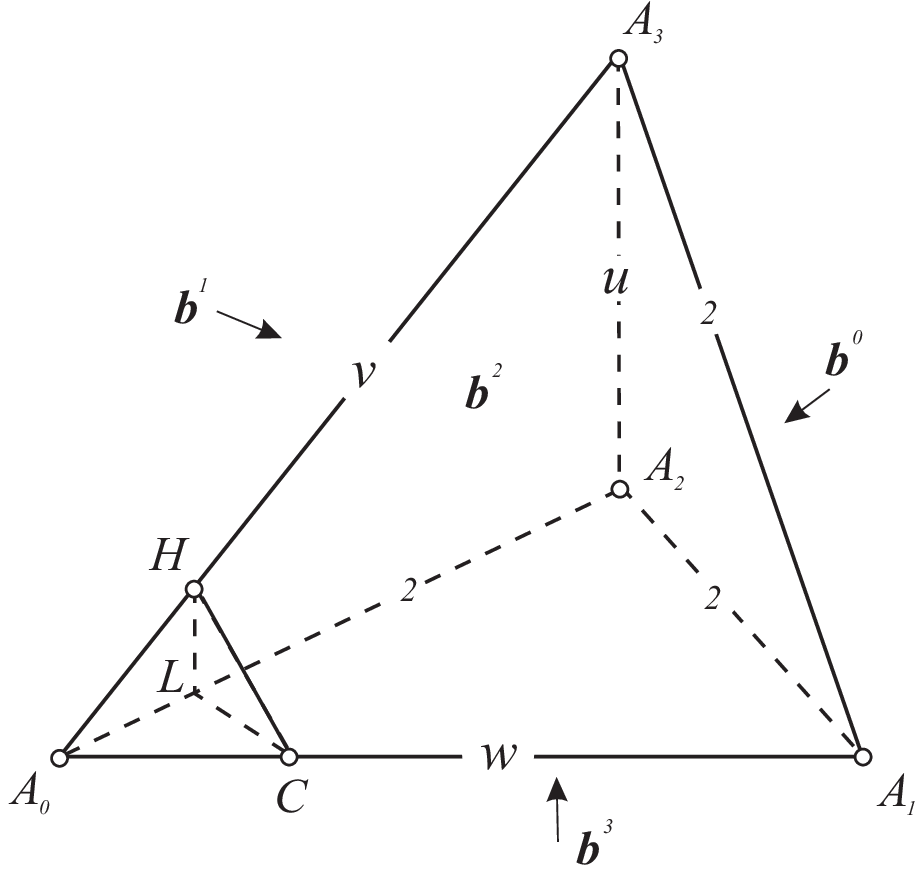} \includegraphics[width=60mm]{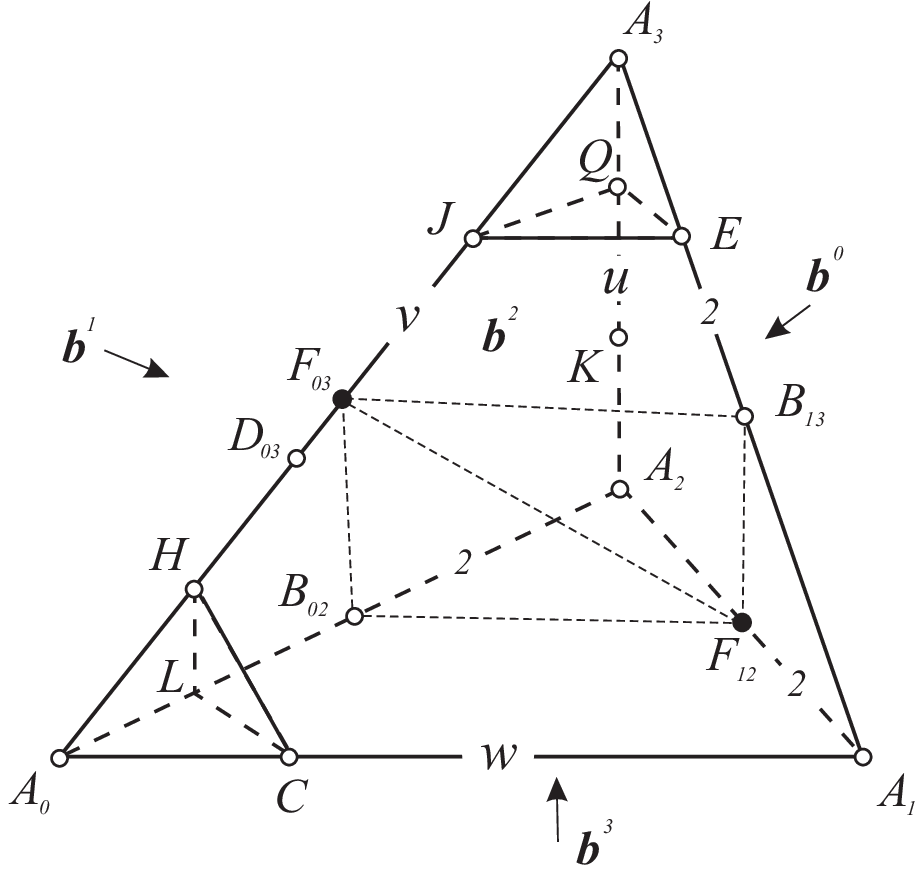}
\caption{Simple and double truncated complete orthoschemes }
\label{Fig1}
\end{figure}
\begin{lem}
Let $A_3$ be an outer principal vertex of the orthoscheme $W_{uvw}$ and let $a_3(\mbox{\boldmath$a$}_{3})=JEQ$ be its polar plane where
$J=a_3 \cap A_3A_0$, $E=a_3 \cap A_3A_1$, $Q=a_3 \cap A_3A_2$ whose vectors are the following:
\begin{equation}
\begin{gathered}
J(\bj)=a_3 \cap A_3A_0; ~ \bj=\ba_0-\frac{a_{03}}{a_{33}} \ba_3, \ \mathrm{with} \\
\langle \bj,\bj \rangle =\frac{(a_{00}a_{33}-a_{03}^2)}{a_{33}}=\langle \bj,\ba_0 \rangle
=\frac{\sin^2\frac{\pi}{v}}{\sin^2\frac{\pi}{u}-\cos^2\frac{\pi}{v}}=\frac{\sin^2\frac{\pi}{v}}{Ba_{33}} \\
E(\mathbf{e})=a_3 \cap A_3A_1; ~ \mathbf{e}=\ba_1-\frac{a_{13}}{a_{33}} \ba_3, \ \mathrm{with} \\
\langle \be,\be \rangle =\frac{(a_{11}a_{33}-a_{13}^2)}{a_{33}}=\langle \be,\ba_1 \rangle
=\frac{1}{\sin^2\frac{\pi}{u}-\cos^2\frac{\pi}{v}}=\frac{1}{Ba_{33}} \\
Q(\bh)=a_3 \cap A_3A_2; ~ \bq=\ba_2-\frac{a_{23}}{a_{33}} \ba_3, \ \mathrm{with} \\
\langle \bq,\bq \rangle =\frac{(a_{22}a_{33}-a_{23}^2)}{a_{33}}=\langle \bq,\ba_2 \rangle
=\frac{\sin^2\frac{\pi}{u}}{\sin^2\frac{\pi}{u}-\cos^2\frac{\pi}{v}}=\frac{a_{22}}{a_{33}}.
\end{gathered} \tag{3.2}
\end{equation}
\end{lem}

\begin{lem}
The midpoint $K(\mathbf{k})$ (see Fig.~1) of $A_2Q$ can be determined by the following vector:
\begin{equation}
\bk=\frac{\ba_2}{\sqrt{-a_{22}}}+\frac{\bq}{\sqrt{-\langle \bq,\bq \rangle}}~ \mathrm{with} ~
\langle \bk,\bk \rangle= -2\Big(1+\sqrt{\frac{1}{a_{33}}}\Big). \tag{3.3}
\end{equation}
\end{lem}
\begin{figure}[htb]
\centering
\includegraphics[width=90mm]{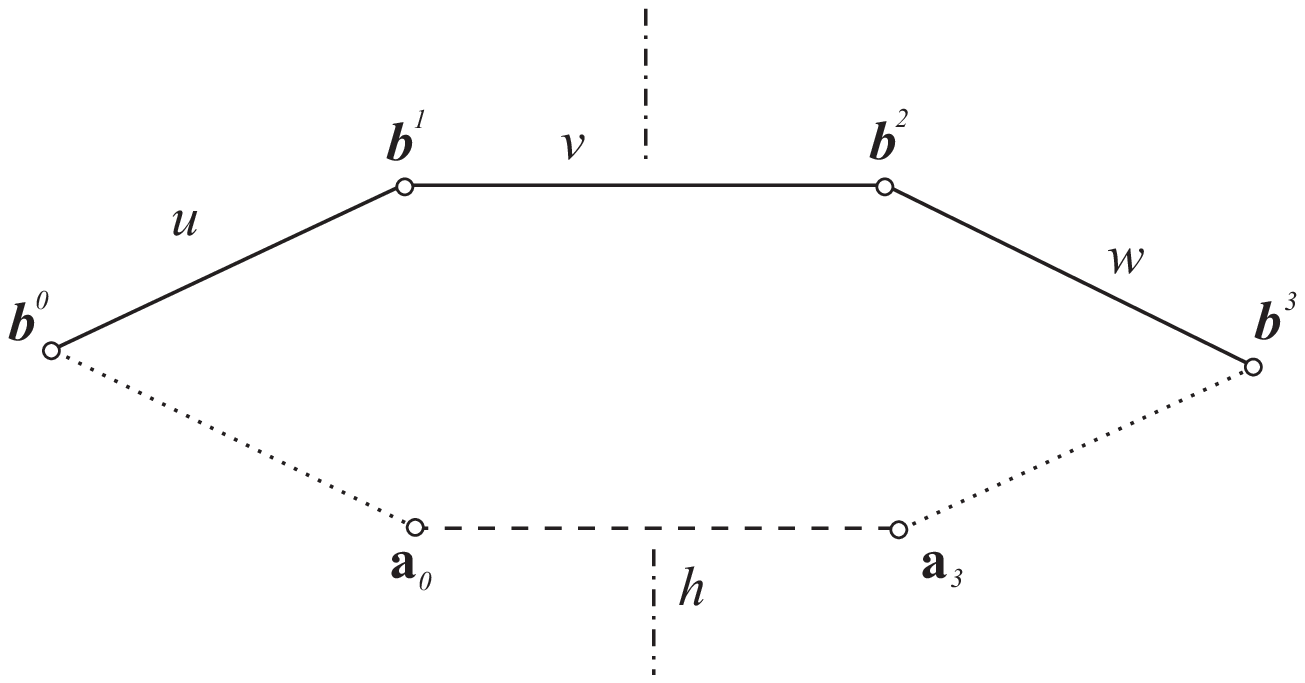}
\caption{}
\label{Fig2}
\end{figure}
Especially, if $u=w$ the midpoints $F_{03}$ of $A_0A_3$ and $F_{12}$ of $A_1A_2$, respectively can play important roles, since
$F_{03}F_{12}$ will be the axis of half turn,
$$
\bh:~ 0 \leftrightarrow 3, 1 \leftrightarrow 2, ~\mathrm{i.e.} ~ A_0 \leftrightarrow A_3, ~
\mbox{\boldmath$b$}^{0} \leftrightarrow \mbox{\boldmath$b$}^{3},  ~ A_1 \leftrightarrow A_2, ~
\mbox{\boldmath$b$}^{1} \leftrightarrow \mbox{\boldmath$b$}^{2}.
$$
(Here $a_{00}=a_{33}$ and $a_{11}=a_{22}$ hold, of course.)
\begin{lem}
The midpoints $F_{03}(\mathbf{f}_{03})$ of $A_3A-0$ and  $F_{12}(\mathbf{f}_{12})$ of $A_1 A_2$ (see Fig.~1) can be
determined by the following vectors:
\begin{equation}
\begin{gathered}
\mathbf{f}_{03}=\ba_0+\ba_3, ~ \langle \mathbf{f}_{03},\mathbf{f}_{03} \rangle= 2(a_{00}+a_{03})<0, \\
\mathbf{f}_{12}=\ba_1+\ba_2, ~ \langle \mathbf{f}_{12},\mathbf{f}_{12} \rangle= 2(a_{11}+a_{12})<0,
\end{gathered} \notag
\end{equation}
independently of $A_0$ and $A_3$ both are either proper, boundary or outer points.
\end{lem}
\begin{figure}[htb]
\centering
\includegraphics[width=80mm]{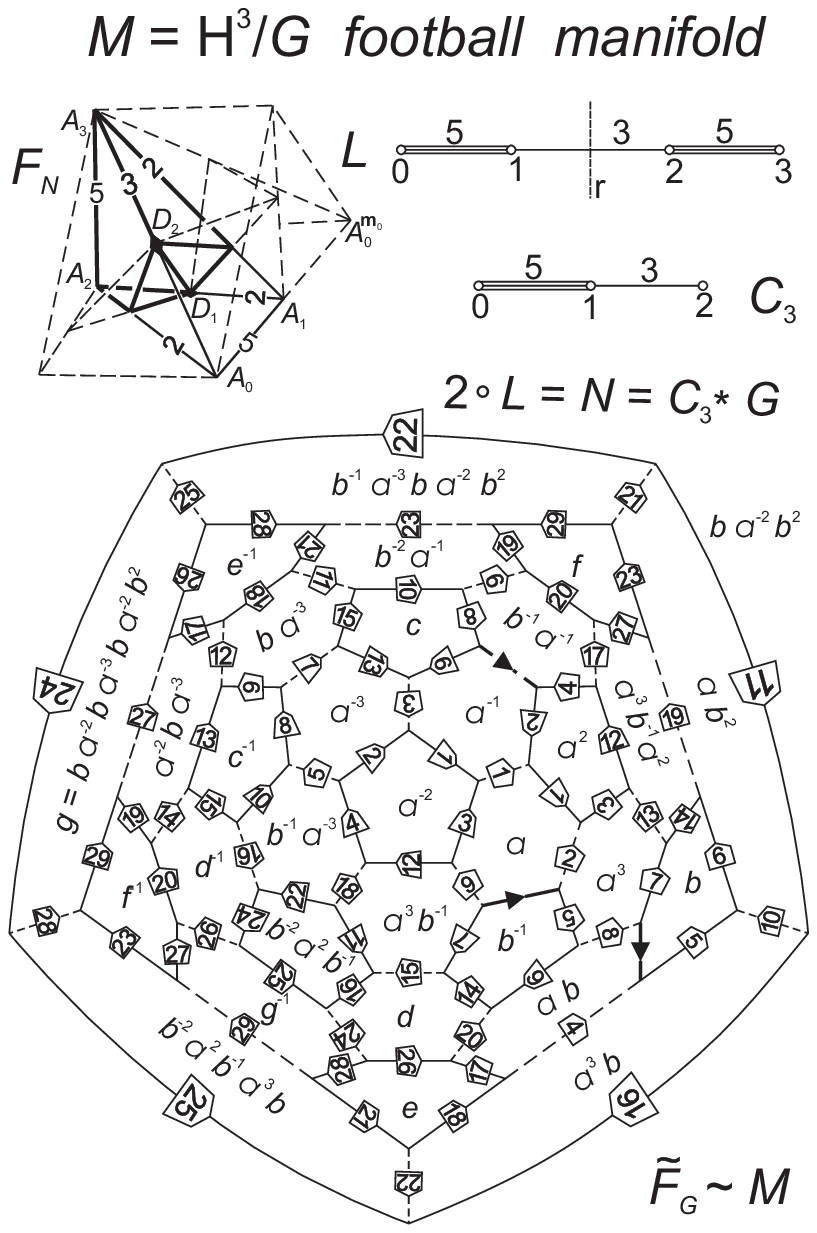}
\caption{}
\label{Fig3}
\end{figure}
The perpendicular foot point $Y(\mathbf{y})$, dropped onto a plane $(\mbox{\boldmath$u$})$ from point $X(\mathbf{x})$ is given by
\begin{equation}
\mathbf{y} = \mathbf{x} - \frac{\langle \mathbf{x}, \mathbf{u} \rangle}{\langle \mathbf{u}, \mathbf{u} \rangle} \mathbf{u}, \tag{3.4}
\end{equation}
where $(\mathbf{u})$ is the pole of the plane $(\mbox{\boldmath$u$})$.
In the considered cases the planes are the faces of the orthoscheme $W_{uvw}$, therefore the poles
$\mbox{\boldmath$b$}^{i}_*=\mathbf{b}^i=b^{ij}\ba_j$ play an important role in the computations that are derived by the matrix (2.1):
\begin{equation}
\begin{gathered}
\mathbf{b}^0=\ba_0-\cos\frac{\pi}{u}\ba_1, \
\mathbf{b}^1=-\cos\frac{\pi}{u}\ba_0+\ba_1-\cos\frac{\pi}{v}\ba_2, \\
\mathbf{b}^2=-\cos\frac{\pi}{v}\ba_1+\ba_2-\cos\frac{\pi}{w}\ba_3, \
\mathbf{b}^3=-\cos\frac{\pi}{w}\ba_2+\ba_3.
\end{gathered} \tag{3.5}
\end{equation}
E.g. the foot point $F_{12}^2(\mathbf{f}_{12}^2)$ of $\mathbf{f}_{12}=\ba_1+\ba_2$ onto the plane $\mbox{\boldmath$b$}^2$ (see Fig.~1) is
$$
\mathbf{f}_{12}^2 = \mathbf{f}_{12}+c\mathbf{b}^2=(1+\cos\frac{\pi}{v})\ba_1+\cos\frac{\pi}{w}\ba_3, ~ \mathrm{by} ~
\mathbf{f}_{12}^2 \mbox{\boldmath$b$}^2=0,
$$
thus for coefficient $c=-1$. We similarly get the other foot points on the simplex side planes $\mbox{\boldmath$b$}^i$ $i \in \{0,1,2,3\}$.

Dividing points of edges of $W_{uvw}$ (or of its half domain) e.g. $D_{ij} \in A_iA_j$ ($i<j\in \{0,1,2,3\}$) can also have special roles with more
cumbersome calculations: $\mathbf{d}_{ij}=\ba_i+d_{ij}\ba_j$ will be interesting for appropriate coeffitient $d_{ij}$. The at least
two locally minimal face distances from a $D_{ij}$ occur for packing. Or at least two locally maximal vertex distances from a $D_{ij}$
occur for covering. Maybe, we obtain different $D_{ij}$'s for locally maximal packing and locally minimal covering, respectively.

We leave out these lengthy discussions here. Only
$F_{03}$, $F_{12}$ and $K$ will be mentioned in Sect.~4 at cases c, d and h, respectively.

Face points $F^i$ on $b^i$ or $F_i$ on $a_i$ can play similar roles (for at least three locally minimal face distances and at least
three locally maximal vertex distances, respectively) for searching locally extremal densities, but we do
not discuss details here, because they will not be relevant for the optimal
packing and covering, respectively.

Extremal "incentre" $I$ in the body of $W_{uvw}$ for packing and "circumcentre" $G$ in $W_{uvw}$ for covering
(for at least four locally minimal face distances and at least
four locally maximal vertex distances, respectively) will not be discussed as well, because they provide only local extrema
(far from our $(5,3,5)$ absolute optima (Table 1.s.i.a)), after lengthy but straightforward discussions.
\section{Ball packings and coverings}
For any fixed $u,v,w$, satisfying the hyperbolity inequlities, we consider the competitive ball centres, as follow,
for maximal packing density and minimal covering density. We indicate the centres where they occur, respectively.
For the ball centre $O$ we have roughly $8$ competitive cases (Fig.~1):
\begin{enumerate}
\item[a.] $A_3$ is the ball centre, $|Stab_{A_3}\bG|=\frac{8uv}{4-(u-2)(v-2)}$,
\item[b.] $A_2$ is the ball centre, $|Stab_{A_2}\bG|=4u$,
\item[c.] $F_{03}$ is the ball centre, $|Stab_{F_{03}}\bG|=4v, ~ \mathrm{if} ~ u=w$,
\item[d.] $F_{12}$ is the ball centre $|Stab_{F_{12}}\bG|=8, ~ \mathrm{if} ~ u=w$,
\item[e.] $Q$ is the ball centre, $|Stab_{Q}\bG|=4u$,
\item[f.] $J$ is the ball centre, $|Stab_{J}\bG|=4v$,
\item[g.] $E$ is the ball centre, $|Stab_{E}\bG|=8$,
\item[h.] the midpoint $K$ of $A_2Q$ is the ball centre, $|Stab_{K}\bG|=2u$.
\end{enumerate}
$K$ gives (only) estimates for the points of $A_2Q$, in general.
Competitive packing ball radii $r$ for $r^{opt}$ with maximal density
\begin{equation}
\delta^{opt}=\frac{\mathrm{Vol}(B(r^{opt}))}{|Stab_{O}\bG|\mathrm{Vol}(W_{uvw})}~ \mathrm{or} ~
\delta^{opt}=\frac{\mathrm{Vol}(B(r^{opt}))}{\frac{1}{2}|Stab_{O}\bG|\mathrm{Vol}(W_{uvw})} ~ \mathrm{if} ~ u=w. \tag{4.1}
\end{equation}
Competive covering ball radii $R$ for $R^{opt}$ with minimal density:
\begin{equation}
\Delta^{opt}=\frac{\mathrm{Vol}(B(R^{opt}))}{|Stab_{O}\bG|\mathrm{Vol}(W_{uvw})}~ \mathrm{or} ~
\Delta^{opt}=\frac{\mathrm{Vol}(B(R^{opt}))}{\frac{1}{2}|Stab_{O}\bG|\mathrm{Vol}(W_{uvw})} ~ \mathrm{if} ~ u=w. \tag{4.2}
\end{equation}
\subsection{Case 1.i.a}
In this case $A_3$ is proper $O=A_3$ is the ball centre. $A_O$ is either proper or boundary point. In the latter case covering is not defined.
The optimal radii can be derived by the {\it Beltrami-Cayley-Klein} model (see Section 2,~3):
\begin{equation}
\begin{gathered}
r=A_3A_2;~ \cosh{(A_3A_2)}=\frac{-a_{23}}{\sqrt{a_{22}a_{33}}}=
\frac{\sin\frac{\pi}{u}\cos\frac{\pi}{w}}{\sqrt{\sin^2\frac{\pi}{u}-\cos^2\frac{\pi}{v}}}=\sqrt{1-\frac{1}{a_{33}}}, \\
R=A_3A_0;~ \cosh{(A_3A_0)}=\frac{-a_{03}}{\sqrt{a_{00}a_{33}}}=\sqrt{1-\frac{\sin^2\frac{\pi}{v}}{Ba_{00}a_{33}}},
\end{gathered} \tag{4.3}
\end{equation}

The densities of the ball packings and coverings for given parameters $u,v,w$ can be computed by the formulas (4.1) and (4.2), also later on.
\medbreak
\scriptsize
\centerline{\vbox{
\halign{\strut\vrule~\hfil $#$ \hfil~\vrule
&\quad \hfil $#$ \hfil~\vrule
&\quad \hfil $#$ \hfil\quad\vrule
&\quad \hfil $#$ \hfil\quad\vrule
&\quad \hfil $#$ \hfil\quad\vrule
\cr
\noalign{\hrule}
\noalign{\vskip2pt}
\multispan5{\strut\vrule\hfill\bf Table 1.i.a, Packing, $O=A_3$, $\frac{1}{u}+\frac{1}{v} > \frac{1}{2},~\frac{1}{v}+\frac{1}{w} \ge \frac{1}{2}$  \hfill\vrule}%
\cr
\noalign{\vskip2pt}
\noalign{\hrule}
\noalign{\vskip2pt}
\noalign{\hrule}
(u,~v,~w) & r^{opt} & Vol(W_{uvw}) & Vol(B(r^{opt})) & \delta^{opt} \cr
\noalign{\hrule}
(3,3,6) & 0.34657 & 0.04229 & 0.17861 & 0.17598 \cr
\noalign{\hrule}
(4,3,5) & 0.53064 & 0.03589 & 0.66207 & 0.38437 \cr
\noalign{\hrule}
(4,3,6) & 0.65847 & 0.10572 & 1.30405 & 0.25697 \cr
\noalign{\hrule}
(5,3,4) & 0.80846 & 0.03589 & 2.52145 & 0.58554 \cr
\noalign{\hrule}
(5,3,5) & 0.99639 & 0.09333 & 5.04848 & 0.45080 \cr
\noalign{\hrule}
(5,3,6) & 1.08394 & 0.17150 & 6.73795 & 0.32740 \cr
\noalign{\hrule}
(3,4,4) & 0.65848 & 0.07633 & 1.30405 & 0.35592 \cr
\noalign{\hrule}
(3,5,3) & 0.86830 & 0.03905 & 3.18663 & {\mathbf{0.68003}} \cr
\noalign{\hrule}}}}
\smallbreak
\scriptsize
\centerline{\vbox{
\halign{\strut\vrule~\hfil $#$ \hfil~\vrule
&\quad \hfil $#$ \hfil~\vrule
&\quad \hfil $#$ \hfil\quad\vrule
&\quad \hfil $#$ \hfil\quad\vrule
&\quad \hfil $#$ \hfil\quad\vrule
\cr
\noalign{\hrule}
\noalign{\vskip2pt}
\multispan5{\strut\vrule\hfill\bf Table 1.i.a, Covering, $O=A_3$, $\frac{1}{u}+\frac{1}{v} \ge \frac{1}{2},~\frac{1}{v}+\frac{1}{w} > \frac{1}{2}$  \hfill\vrule}%
\cr
\noalign{\vskip2pt}
\noalign{\hrule}
\noalign{\vskip2pt}
\noalign{\hrule}
(u,~v,~w) & R^{opt} & Vol(W_{uvw}) & Vol(B(r^{opt})) & \Delta^{opt} \cr
\noalign{\hrule}
(4,3,5) & 1.22646 & 0.03589 & 10.41484& 6.04641 \cr
\noalign{\hrule}
(5,3,4) & 1.22646 & 0.03589 & 10.41484 & {\mathbf{2.41856}} \cr
\noalign{\hrule}
(5,3,5) & 1.90285 & 0.09333 & 58.62658 & 5.23495 \cr
\noalign{\hrule}
(3,5,3) & 1.38257 & 0.03905 & 16.16044 & 3.44864 \cr
\noalign{\hrule}}}}
\smallbreak
\normalsize

\subsection{Case 1.s.i.a}
$O=A_3$ is the ball centre and
$u=w$. The density is related to the "half orthoscheme".
\begin{equation}
\begin{gathered}
r=\min\{A_3A_2,A_3F_{03}\};~ \mathrm{where} ~ \cosh{(A_3A_2)} ~ \text{is in (4.3)}, \\
\cosh{(A_3F_{03})}=\sqrt{\frac{a_{03}}{2a_{33}}+\frac{1}{2}}=\cosh{\Big(\frac{1}{2} A_3A_0\Big)}, \\
R=A_3F_{12};~ \cosh{(A_3F_{12})}=\frac{-(a_{13}+a_{23})}{\sqrt{2a_{33}(a_{12}+a_{22})}}=\sqrt{\frac{1}{2}+\frac{a_{03}-1}{2a_{33}}},
\end{gathered} \tag{4.4}
\end{equation}
\medbreak
\scriptsize
\centerline{\vbox{
\halign{\strut\vrule~\hfil $#$ \hfil~\vrule
&\quad \hfil $#$ \hfil~\vrule
&\quad \hfil $#$ \hfil\quad\vrule
&\quad \hfil $#$ \hfil\quad\vrule
&\quad \hfil $#$ \hfil\quad\vrule
\cr
\noalign{\hrule}
\noalign{\vskip2pt}
\multispan5{\strut\vrule\hfill\bf Table 1.s.i.a, Packing, $O=A_3$, $\frac{1}{u}+\frac{1}{v} > \frac{1}{2},~\frac{1}{v}+\frac{1}{w} \ge \frac{1}{2}$  \hfill\vrule}%
\cr
\noalign{\vskip2pt}
\noalign{\hrule}
\noalign{\vskip2pt}
\noalign{\hrule}
(u,~v,~w) & r^{opt} & Vol(W_{uvw}) & Vol(B(r^{opt})) & \delta^{opt} \cr
\noalign{\hrule}
(5,3,5) & 0.95142 & 0.09333 & 4.31988 & \mathbf{0.77147} \cr
\noalign{\hrule}
(3,5,3) & 0.69129 & 0.03905 & 1.52220 & 0.64967 \cr
\noalign{\hrule}}}}
\smallbreak
\scriptsize
\centerline{\vbox{
\halign{\strut\vrule~\hfil $#$ \hfil~\vrule
&\quad \hfil $#$ \hfil~\vrule
&\quad \hfil $#$ \hfil\quad\vrule
&\quad \hfil $#$ \hfil\quad\vrule
&\quad \hfil $#$ \hfil\quad\vrule
\cr
\noalign{\hrule}
\noalign{\vskip2pt}
\multispan5{\strut\vrule\hfill\bf Table 1.s.i.a, Covering, $O=A_3$, $\frac{1}{u}+\frac{1}{v} \ge \frac{1}{2},~\frac{1}{v}+\frac{1}{w} > \frac{1}{2}$  \hfill\vrule}%
\cr
\noalign{\vskip2pt}
\noalign{\hrule}
\noalign{\vskip2pt}
\noalign{\hrule}
(u,~v,~w) & R^{opt} & Vol(W_{uvw}) & Vol(B(R^{opt})) & \Delta^{opt} \cr
\noalign{\hrule}
(5,3,5) & 1.12484 & 0.09333 & 7.66539 & {\mathbf{1.36893}} \cr
\noalign{\hrule}
(3,5,3) & 0.89558 & 0.03905 & 3.53002 & 1.50661 \cr
\noalign{\hrule}}}}
\smallbreak
\normalsize
\subsection{Case 1.i.b}
$O=A_2$ is the ball centre.
\begin{equation}
\begin{gathered}
r=A_2b^2=:A_2A_2^2;~ \cosh{(A_2b^2)}=\sqrt{1-\frac{1}{a_{22}}} \\
R=\max\{A_2A_3,A_2A_0\};~ \cosh{(A_2A_3)}=\sqrt{1-\frac{1}{a_{33}}},\\
\cosh{(A_2A_0)}=\frac{-a_{02}}{\sqrt{a_{00}a_{22}}}=\frac{\cot\frac{\pi}{u}\cos\frac{\pi}{v}}{\sqrt{\sin^2\frac{\pi}{w}-\cos^2\frac{\pi}{v}}}
=\sqrt{1-\frac{1}{Ba_{00}a_{22}}}.
\end{gathered} \tag{4.5}
\end{equation}

\medbreak
\scriptsize
\centerline{\vbox{
\halign{\strut\vrule~\hfil $#$ \hfil~\vrule
&\quad \hfil $#$ \hfil~\vrule
&\quad \hfil $#$ \hfil\quad\vrule
&\quad \hfil $#$ \hfil\quad\vrule
&\quad \hfil $#$ \hfil\quad\vrule
\cr
\noalign{\hrule}
\noalign{\vskip2pt}
\multispan5{\strut\vrule\hfill\bf Table 1.i.b, Packing, $O=A_2$, $\frac{1}{u}+\frac{1}{v} \ge \frac{1}{2},~\frac{1}{v}+\frac{1}{w} \ge \frac{1}{2}$  \hfill\vrule}%
\cr
\noalign{\vskip2pt}
\noalign{\hrule}
\noalign{\vskip2pt}
\noalign{\hrule}
(u,~v,~w) & r^{opt} & Vol(W_{uvw}) & Vol(B(r^{opt})) & \delta^{opt} \cr
\noalign{\hrule}
(4,3,5) & 0.38360 & 0.03589 & 0.24350 & 0.42409 \cr
\noalign{\hrule}
(5,3,4) & 0.45682 & 0.03589 & 0.41631 & {\mathbf{0.58007}} \cr
\noalign{\hrule}
(5,3,5) & 0.58157 & 0.09333 & 0.88150 & 0.47227 \cr
\noalign{\hrule}
(4,4,3) & 0.48121 & 0.07633 & 0.48886 & 0.40029 \cr
\noalign{\hrule}
(3,5,3) & 0.34346 & 0.03905 & 0.17377 & 0.37082 \cr
\noalign{\hrule}}}}
\smallbreak
\scriptsize
\centerline{\vbox{
\halign{\strut\vrule~\hfil $#$ \hfil~\vrule
&\quad \hfil $#$ \hfil~\vrule
&\quad \hfil $#$ \hfil\quad\vrule
&\quad \hfil $#$ \hfil\quad\vrule
&\quad \hfil $#$ \hfil\quad\vrule
\cr
\noalign{\hrule}
\noalign{\vskip2pt}
\multispan5{\strut\vrule\hfill\bf Table 1.i.b, Covering, $O=A_2$, $\frac{1}{u}+\frac{1}{v} \ge \frac{1}{2},~\frac{1}{v}+\frac{1}{w} > \frac{1}{2}$  \hfill\vrule}%
\cr
\noalign{\vskip2pt}
\noalign{\hrule}
\noalign{\vskip2pt}
\noalign{\hrule}
(u,~v,~w) & R^{opt} & Vol(W_{uvw}) & Vol(B(r^{opt})) & \Delta^{opt} \cr
\noalign{\hrule}
(5,3,4) & 0.84248 & 0.03589 & 2.88536 & {\mathbf{4.02028}} \cr
\noalign{\hrule}}}}
\smallbreak
\normalsize
\subsection{Case 1.s.i.b}
$O=A_2$ is the ball centre.
$u=w$. The density is related to the "half orthoscheme".
\begin{equation}
\begin{gathered}
r=\min\{A_2b^2,A_2F_{12}\};~ \mathrm{where} ~ \cosh{(A_2b^2)}~ \text{in (4.5)}, \\
\cosh{(A_2F_{12})}=\sqrt{\frac{a_{12}}{2a_{22}}+\frac{1}{2}}=\cosh{\Big(\frac{1}{2} A_2A_1\Big)}=\sqrt{\frac{1}{2}+\frac{\cos\frac{\pi}{v}}{2\sin^2\frac{\pi}{u}}}, \\
R=\max\{A_2 A_3,A_2F_{03}\};~
\cosh{(A_2A_3)}=\frac{-a_{23}}{\sqrt{a_{22}a_{33}}}= \\
\cosh{(A_2F_{03})}=\frac{-(a_{02}+a_{23})}{\sqrt{2a_{22}(a_{33}+a_{03})}}=\cot\frac{\pi}{u}\sqrt{\frac{\cos\frac{\pi}{v}+\sin^2\frac{\pi}{u}}
{2(1-\cos\frac{\pi}{v})}}.
\end{gathered} \tag{4.6}
\end{equation}

\medbreak
\scriptsize
\centerline{\vbox{
\halign{\strut\vrule~\hfil $#$ \hfil~\vrule
&\quad \hfil $#$ \hfil~\vrule
&\quad \hfil $#$ \hfil\quad\vrule
&\quad \hfil $#$ \hfil\quad\vrule
&\quad \hfil $#$ \hfil\quad\vrule
\cr
\noalign{\hrule}
\noalign{\vskip2pt}
\multispan5{\strut\vrule\hfill\bf Table 1.s.i.b, Packing, $O=A_2$, $\frac{1}{u}+\frac{1}{v} \ge \frac{1}{2},~\frac{1}{v}+\frac{1}{w} \ge \frac{1}{2}$  \hfill\vrule}%
\cr
\noalign{\vskip2pt}
\noalign{\hrule}
\noalign{\vskip2pt}
\noalign{\hrule}
(u,~v,~w) & r^{opt} & Vol(W_{uvw}) & Vol(B(r^{opt})) & \delta^{opt} \cr
\noalign{\hrule}
(5,3,5) & 0.45682 & 0.09333 & 0.41631 & {\mathbf{0.44609}} \cr
\noalign{\hrule}}}}
\smallbreak
\scriptsize
\centerline{\vbox{
\halign{\strut\vrule~\hfil $#$ \hfil~\vrule
&\quad \hfil $#$ \hfil~\vrule
&\quad \hfil $#$ \hfil\quad\vrule
&\quad \hfil $#$ \hfil\quad\vrule
&\quad \hfil $#$ \hfil\quad\vrule
\cr
\noalign{\hrule}
\noalign{\vskip2pt}
\multispan5{\strut\vrule\hfill\bf Table 1.s.i.b, Covering, $O=A_2$, $\frac{1}{u}+\frac{1}{v} \ge \frac{1}{2},~\frac{1}{v}+\frac{1}{w} \ge \frac{1}{2}$  \hfill\vrule}%
\cr
\noalign{\vskip2pt}
\noalign{\hrule}
\noalign{\vskip2pt}
\noalign{\hrule}
(u,~v,~w) & R^{opt} & Vol(W_{uvw}) & Vol(B(R^{opt})) & \Delta^{opt} \cr
\noalign{\hrule}
(5,3,5) & 0.99639 & 0.09333 & 5.04848 & {\mathbf{5.40954}} \cr
\noalign{\hrule}}}}
\smallbreak
\normalsize
\subsection{Case 1.s.i.c}
$O=F_{03}$ is the ball centre,
$u=w$.
\begin{equation}
\begin{gathered}
r=F_{03} b^0 = F_{03} b^3;~ \mathrm{where} ~ \cosh{(F_{03} b^0)}=\sqrt{1+\frac{\cos\frac{\pi}{v}-\sin^2\frac{\pi}{u}}{2(1-\cos\frac{\pi}{v})}}, \\
R=\max\{F_{03} A_3, F_{03}A_2\}=F_{03} A_3;~
\cosh{(F_{03} A_3)}=\sqrt{\frac{a_{03}}{2a_{33}}+\frac{1}{2}}.
\end{gathered} \tag{4.7}
\end{equation}

\medbreak
\scriptsize
\centerline{\vbox{
\halign{\strut\vrule~\hfil $#$ \hfil~\vrule
&\quad \hfil $#$ \hfil~\vrule
&\quad \hfil $#$ \hfil\quad\vrule
&\quad \hfil $#$ \hfil\quad\vrule
&\quad \hfil $#$ \hfil\quad\vrule
\cr
\noalign{\hrule}
\noalign{\vskip2pt}
\multispan5{\strut\vrule\hfill\bf Table 1.s.i.c, Packing, $O=F_{03}$, $u=w$,~ $\frac{1}{u}+\frac{1}{v} \ge \frac{1}{2},~\frac{1}{v}+\frac{1}{w} \ge \frac{1}{2}$  \hfill\vrule}%
\cr
\noalign{\vskip2pt}
\noalign{\hrule}
\noalign{\vskip2pt}
\noalign{\hrule}
(u,~v,~w) & r^{opt} & Vol(W_{uvw}) & Vol(B(r^{opt})) & \delta^{opt} \cr
\noalign{\hrule}
(4,4,4) & 0.56419 & 0.22899 & 0.80163 & 0.43759 \cr
\noalign{\hrule}
(3,5,3) & 0.38360 & 0.03905 & 0.24350 & {\mathbf{0.62355}} \cr
\noalign{\hrule}
(3,6,3) & 0.61795 & 0.16916 & 1.06673 & 0.52551 \cr
\noalign{\hrule}}}}
\smallbreak
\scriptsize
\centerline{\vbox{
\halign{\strut\vrule~\hfil $#$ \hfil~\vrule
&\quad \hfil $#$ \hfil~\vrule
&\quad \hfil $#$ \hfil\quad\vrule
&\quad \hfil $#$ \hfil\quad\vrule
&\quad \hfil $#$ \hfil\quad\vrule
\cr
\noalign{\hrule}
\noalign{\vskip2pt}
\multispan5{\strut\vrule\hfill\bf Table 1.s.i.c, Covering, $O=F_{03}$, $u=w$,~  $\frac{1}{u}+\frac{1}{v} \ge \frac{1}{2},~\frac{1}{v}+\frac{1}{w} > \frac{1}{2}$  \hfill\vrule}%
\cr
\noalign{\vskip2pt}
\noalign{\hrule}
\noalign{\vskip2pt}
\noalign{\hrule}
(u,~v,~w) & R^{opt} & Vol(W_{uvw}) & Vol(B(R^{opt})) & \Delta^{opt} \cr
\noalign{\hrule}
(3,5,3) & 0.69129 & 0.03905 & 1.52220 & {\mathbf{3.89804}} \cr
\noalign{\hrule}}}}
\smallbreak
\normalsize
\subsection{Case 1.s.i.d}
$O=F_{12}$ is the ball centre and
$u=w$.
\begin{equation}
\begin{gathered}
r=F_{12}b^1=F_{12}b^2;~ \mathrm{where} ~ \cosh{(F_{12}b^2)}=\sqrt{1+\frac{\cos\frac{\pi}{v}-\sin^2\frac{\pi}{u}}{2}}, \\
R=A_3F_{12};~ \cosh{(A_3F_{12})}=\frac{-(a_{13}+a_{23})}{\sqrt{2a_{33}(a_{12}+a_{22})}}=\sqrt{\frac{1}{2}+\frac{a_{03}-1}{2a_{33}}}.
\end{gathered} \tag{4.8}
\end{equation}
We do not obtain relevant arrangements for optima.
\normalsize
\subsection{Case 1.ii.a}
In this case $A_3$ is proper $O=A_3$ is the ball centre. $A_O$ is outer, $a_0=CLH$ (see Fig.~1) is its polar plane.
\begin{equation}
\begin{gathered}
r=\min\{A_3A_2,A_3H\};~ \cosh{(A_3A_2)}~ \text{is in (4.5)}, \\
\cosh{(A_3H)}=\sqrt{1-\frac{a_{03}^2}{a_{00}a_{33}}}=\frac{\sin\frac{\pi}{v}}{\sqrt{Ba_{00}a_{33}}}, \\
R=A_3C;~ \cosh{(A_3C)}={\frac{a_{01} a_{03} - a_{00} a_{13}}{\sqrt{a_{00}a_{33}(a_{11}a_{00}-a_{01}^2)}}}
=\frac{\cot{\frac{\pi}{w}}\cos \frac{\pi}{v}}{\sqrt{a_{00}a_{33}B}}= \\
=\sqrt{1-\frac{1}{Ba_{11}a_{33}}-\frac{a_{03}^2}{a_{00} a_{33}}}.
\end{gathered} \tag{4.9}
\end{equation}

\medbreak
{\scriptsize{
\centerline{\vbox{
\halign{\strut\vrule~\hfil $#$ \hfil~\vrule
&\quad \hfil $#$ \hfil~\vrule
&\quad \hfil $#$ \hfil\quad\vrule
&\quad \hfil $#$ \hfil\quad\vrule
&\quad \hfil $#$ \hfil\quad\vrule
\cr
\noalign{\hrule}
\noalign{\vskip2pt}
\multispan5{\strut\vrule\hfill\bf Table 1.ii.a, Packing, $O=A_3$, $\frac{1}{u}+\frac{1}{v}>\frac{1}{2},~\frac{1}{v}+\frac{1}{w}<\frac{1}{2}$  \hfill\vrule}%
\cr
\noalign{\vskip2pt}
\noalign{\hrule}
\noalign{\vskip2pt}
\noalign{\hrule}
(u,~v,~w) & r^{opt} & Vol(W_{uvw}) & Vol(B(r^{opt})) & \delta^{opt} \cr
\noalign{\hrule}
(3,3, w \to \infty) & \begin{gathered} \mathrm{arsh}(\sqrt{2}/{2}) \approx \\ \approx 0.65848 \end{gathered}&
0.15266 & 1.30405 & \mathbf{0.35592} \cr
\noalign{\hrule}
(3,4,6) & 0.96242 & 0.19616 & 4.49014 & \mathbf{0.47687} \cr
\noalign{\hrule}
(3,4, w \to \infty) & \begin{gathered} \log(1+\sqrt{2}) \approx \\ \approx 0.88137 \end{gathered}& 0.25096 & 3.34793 & 0.27793 \cr
\noalign{\hrule}
(3,5,4) & 1.30631 & 0.21299 & 13.09457 & 0.51233 \cr
\noalign{\hrule}
(3,5,5) & 1.40036 & 0.26320 & 16.95557 & \mathbf{0.53684} \cr
\noalign{\hrule}
(3,5, w \to \infty) & \begin{gathered} {-\mathrm{arsh}}\Bigg(\frac{-1}{\sqrt{2-2\cos({\pi}/{5})}} \Bigg) \approx \\ \approx 1.25850 \end{gathered}
& 0.33233 & 11.43025 & 0.28662 \cr
\noalign{\hrule}
(4,3, w \to \infty) & \begin{gathered} \log(1+\sqrt{2}) \approx \\ \approx 0.88137 \end{gathered}
& 0.25096 & 3.34793 & \mathbf{0.27793} \cr
\noalign{\hrule}
(5,3, w \to \infty) & \begin{gathered} {-\mathrm{arsh}}\Bigg(\frac{-1}{\sqrt{2-2\cos({\pi}/{5})}} \Bigg) \approx \\ \approx 1.25850 \end{gathered}
& 0.33233 & 11.43025 & \mathbf{0.28662} \cr
\noalign{\hrule}}}}}}
\smallbreak
{\scriptsize{
\centerline{\vbox{
\halign{\strut\vrule~\hfil $#$ \hfil~\vrule
&\quad \hfil $#$ \hfil~\vrule
&\quad \hfil $#$ \hfil\quad\vrule
&\quad \hfil $#$ \hfil\quad\vrule
&\quad \hfil $#$ \hfil\quad\vrule
\cr
\noalign{\hrule}
\noalign{\vskip2pt}
\multispan5{\strut\vrule\hfill\bf Table 1.ii.a, Covering, $O=A_3$, $\frac{1}{u}+\frac{1}{v}>\frac{1}{2},~\frac{1}{v}+\frac{1}{w}<\frac{1}{2}$  \hfill\vrule}%
\cr
\noalign{\vskip2pt}
\noalign{\hrule}
\noalign{\vskip2pt}
\noalign{\hrule}
(u,~v,~w) & R^{opt} & Vol(W_{uvw}) & Vol(B(R^{opt})) & \Delta^{opt} \cr
\noalign{\hrule}
(3,3,8) & 1.30889 & 0.10721 & 13.18949 & {\mathbf{5.12591}} \cr
\noalign{\hrule}
(3,4,5) & 1.53591 & 0.16596 & 24.17649 & {\mathbf{3.03486}} \cr
\noalign{\hrule}
(3,5,4) & 1.93116 & 0.21299 & 62.56492 & {\mathbf{2.44790}} \cr
\noalign{\hrule}
(3,5,5) & 2.08287 & 0.26320 & 88.11084 & 2.78974 \cr
\noalign{\hrule}
(4,3,8) & 1.81579 & 0.18790 & 47.88239 & {\mathbf{5.30903}} \cr
\noalign{\hrule}
(5,3,8) & 2.37474 & 0.26094& 166.52911 & {\mathbf{5.31816}} \cr
\noalign{\hrule}}}}}}
\smallbreak

\normalsize
\subsection{Case 1.ii.b}
$O=A_2$ is the ball centre.
\begin{equation}
\begin{gathered}
r=\min\{A_2b^2,A_2L\};~ \mathrm{where} ~ \cosh{(A_2b^2)}=\sqrt{1-\frac{1}{a_{22}}}, \\
A_2a_0=A_2L;~ \cosh{(A_2L)}=\sqrt{1-\frac{a_{02}^2}{a_{00} a_{22}}}
=\frac{1}{\sin{u}}\sqrt{\frac{B}{\sin^2\frac{\pi}{w}-\cos^2\frac{\pi}{v}}}, \\
R=\max\{A_2C,A_2A_3,A_2H\};~ \cosh{(A_2C)}={\frac{a_{01} a_{02} - a_{12} a_{00}}{\sqrt{a_{00}a_{22}(a_{11}a_{00}-a_{01}^2)}}}
=\\ =\frac{\cos{\frac{\pi}{v}}}{\sin \frac{\pi}{u}\sin \frac{\pi}{w}\sqrt{a_{00}}}=
\sqrt{1-\frac{1}{Ba_{11}a_{22}}-\frac{a_{02}^2a_{12}}{a_{00} a_{22}}}, \\
\cosh{(A_2H)}={\frac{a_{02} a_{03} - a_{23} a_{00}}{\sqrt{a_{00}a_{22}(a_{33}a_{00}-a_{03}^2)}}}
=\frac{\cot{\frac{\pi}{w}}}{\sin \frac{\pi}{v}\sqrt{a_{00}}}, \\
\cosh{(A_2A_3)} ~ \text{is in (4.5).}
\end{gathered} \tag{4.10}
\end{equation}

\medbreak
{\scriptsize{
\centerline{\vbox{
\halign{\strut\vrule~\hfil $#$ \hfil~\vrule
&\quad \hfil $#$ \hfil~\vrule
&\quad \hfil $#$ \hfil\quad\vrule
&\quad \hfil $#$ \hfil\quad\vrule
&\quad \hfil $#$ \hfil\quad\vrule
\cr
\noalign{\hrule}
\noalign{\vskip2pt}
\multispan5{\strut\vrule\hfill\bf Table 1.ii.b, Packing, $O=A_2$, $\frac{1}{u}+\frac{1}{v} \ge \frac{1}{2},~\frac{1}{v}+\frac{1}{w}<\frac{1}{2}$  \hfill\vrule}%
\cr
\noalign{\vskip2pt}
\noalign{\hrule}
\noalign{\vskip2pt}
\noalign{\hrule}
(u,~v,~w) & r^{opt} & Vol(W_{uvw}) & Vol(B(r^{opt})) & \delta^{opt} \cr
\noalign{\hrule}
(3,4,6) & 0.60745 & 0.19616 & 1.01065 & {0.42934} \cr
\noalign{\hrule}
(3,4,7) & 0.64583 & 0.21218 & 1.22631 & {\mathbf{0.48164}} \cr
\noalign{\hrule}
(3,3, w \to \infty) & \begin{gathered} \mathrm{arch}(2/\sqrt{3}) \approx \\ \approx 0.54931 \end{gathered}&
0.15266 & 0.73740 & \mathbf{0.40253} \cr
\noalign{\hrule}
(3,5,5) & 0.67390 & 0.26320 & 1.40355 & \mathbf{0.44439} \cr
\noalign{\hrule}
(3,5, w \to \infty) & \begin{gathered} {\mathrm{arch}}(2/\sqrt{3}) \approx \\ \approx 0.54931 \end{gathered}
& 0.33233 & 0.73740 & 0.18491 \cr
\noalign{\hrule}
(4,3, w \to \infty) & \begin{gathered} {\mathrm{arch}(\sqrt{6}/2)} \approx \\ \approx 0.65848 \end{gathered}
& 0.25096 & 1.30405 & \mathbf{0.32477} \cr
\noalign{\hrule}
(5,3, w \to \infty) & \begin{gathered} {\mathrm{arch}}\Bigg(\frac{\sqrt{1+4 \sin^2(\pi/5)}}{\sin(\pi/5)} \Bigg) \approx \\ \approx 0.77173 \end{gathered}
& 0.33233 & 2.16804 & \mathbf{0.32619} \cr
\noalign{\hrule}
(3,6,5) & 0.72182 & 0.35992 & 1.74787 & {\mathbf{0.40469}} \cr
\noalign{\hrule}
(3,6, w \to \infty) & \begin{gathered} {\mathrm{arch}}(2/\sqrt{3}) \approx \\ \approx 0.54931 \end{gathered}
& 0.42289 & 0.73740 & 0.14531 \cr
\noalign{\hrule}
(6,3, w \to \infty) & \begin{gathered} {\mathrm{arch}}(\sqrt{2}) \approx \\ \approx 0.88137 \end{gathered}
& 0.42289 & 3.34793 & \mathbf{0.32987} \cr
\noalign{\hrule}
(4,4, w \to \infty) & \begin{gathered} {\mathrm{arch}}(\sqrt{2}) \approx \\ \approx 0.88137 \end{gathered}
& 0.45798 & 3.34793 & \mathbf{0.45689} \cr
\noalign{\hrule}}}}}}
\smallbreak
{\scriptsize{
\centerline{\vbox{
\halign{\strut\vrule~\hfil $#$ \hfil~\vrule
&\quad \hfil $#$ \hfil~\vrule
&\quad \hfil $#$ \hfil\quad\vrule
&\quad \hfil $#$ \hfil\quad\vrule
&\quad \hfil $#$ \hfil\quad\vrule
\cr
\noalign{\hrule}
\noalign{\vskip2pt}
\multispan5{\strut\vrule\hfill\bf Table 1.ii.b, Covering, $O=A_2$, $\frac{1}{u}+\frac{1}{v}>\frac{1}{2},~\frac{1}{v}+\frac{1}{w}<\frac{1}{2}$  \hfill\vrule}%
\cr
\noalign{\vskip2pt}
\noalign{\hrule}
\noalign{\vskip2pt}
\noalign{\hrule}
(u,~v,~w) & R^{opt} & Vol(W_{uvw}) & Vol(B(R^{opt})) & \Delta^{opt} \cr
\noalign{\hrule}
(3,3,8) & 1.16245 & 0.10721 & 8.60501 & {\mathbf{6.68843}} \cr
\noalign{\hrule}
(3,4,5) & 1.14791 & 0.16596 & 8.23135 & {\mathbf{4.13311}} \cr
\noalign{\hrule}
(3,5,4) & 1.30632 & 0.21299 & 13.09457 & {\mathbf{5.12335}} \cr
\noalign{\hrule}
(4,3,9) & 1.55832 & 0.20295 & 25.59284 & {\mathbf{7.88150}} \cr
\noalign{\hrule}
(5,3,9) & 1.80489 & 0.27783 & 46.67018 & {\mathbf{8.39914}} \cr
\noalign{\hrule}}}}}}
\smallbreak
\normalsize
\subsection{Case 2.i.b}
In this case $A_3$ is outer, $a_3=JEQ$ is its polar plane, $A_0$ is proper or boundary point. $O=A_2$ is the ball centre (see Fig.~1).
\begin{equation}
\begin{gathered}
r=\min\{A_2b^2,A_2a_3\};~ \cosh{(A_2b^2)}=\sqrt{1-\frac{1}{a_{22}}},\\
A_2a_3=A_2Q;~ \cosh{(A_2Q)}=\frac{1}{\sqrt{a_{33}}}
=\sqrt{1+\frac{\sin^2\frac{\pi}{u}\cos^2\frac{\pi}{w}}{\cos^2\frac{\pi}{v}-\sin^2\frac{\pi}{u}}}, \\
R=\max\{A_2A_0,A_2J\};~ \cosh{(A_2A_0)}=\frac{-a_{02}}{\sqrt{a_{00}a_{22}}},\\
\cosh{(A_2J)}={\frac{a_{03} a_{23} - a_{02} a_{33}}{\sqrt{a_{33}a_{22}(a_{33}a_{00}-a_{03}^2)}}}
=\frac{\cot{\frac{\pi}{u}}\cot\frac{\pi}{v}}{\sqrt{a_{00}}}.
\end{gathered} \tag{4.11}
\end{equation}
\medbreak
\scriptsize
\centerline{\vbox{
\halign{\strut\vrule~\hfil $#$ \hfil~\vrule
&\quad \hfil $#$ \hfil~\vrule
&\quad \hfil $#$ \hfil\quad\vrule
&\quad \hfil $#$ \hfil\quad\vrule
&\quad \hfil $#$ \hfil\quad\vrule
\cr
\noalign{\hrule}
\noalign{\vskip2pt}
\multispan5{\strut\vrule\hfill\bf Table 2.i.b. Packing, $O=A_2$, $\frac{1}{u}+\frac{1}{v} < \frac{1}{2},~\frac{1}{v}+\frac{1}{w} \ge \frac{1}{2}$  \hfill\vrule}%
\cr
\noalign{\vskip2pt}
\noalign{\hrule}
\noalign{\vskip2pt}
\noalign{\hrule}
(u,~v,~w) & r^{opt} & Vol(W_{uvw}) & Vol(B(r^{opt})) & \delta^{opt} \cr
\noalign{\hrule}
(7,3,3) & 0.70133 & 0.08856 & 1.59395 & \mathbf{0.64279} \cr
\noalign{\hrule}
(7,3,4) & 0.81624 & 0.16297 & 2.60141 & \mathbf{0.57008} \cr
\noalign{\hrule}
(7,3,5) & 0.87514 & 0.23326 & 3.27033 & \mathbf{0.50072} \cr
\noalign{\hrule}
(5,4,3) & 0.69129 & 0.16596 & 1.52220 & \mathbf{0.45859} \cr
\noalign{\hrule}
(4,5,3) & 0.69129 & 0.21299 & 1.52220 & \mathbf{0.44668} \cr
\noalign{\hrule}
(7,3,6) & 0.90817 & 0.31781 & 3.69768 & \mathbf{0.41553} \cr
\noalign{\hrule}
(5,4,4) & 0.86233 & 0.34084 & 3.11499 & \mathbf{0.45696} \cr
\noalign{\hrule}
(4,6,3) & 0.65848 & 0.31717 & 1.30405 & \mathbf{0.25697} \cr
\noalign{\hrule}}}}
\smallbreak
\centerline{\vbox{
\halign{\strut\vrule~\hfil $#$ \hfil~\vrule
&\quad \hfil $#$ \hfil~\vrule
&\quad \hfil $#$ \hfil\quad\vrule
&\quad \hfil $#$ \hfil\quad\vrule
&\quad \hfil $#$ \hfil\quad\vrule
\cr
\noalign{\hrule}
\noalign{\vskip2pt}
\multispan5{\strut\vrule\hfill\bf Table 2.i.b. Covering, $O=A_2$, $\frac{1}{u}+\frac{1}{v}<\frac{1}{2},~\frac{1}{v}+\frac{1}{w} > \frac{1}{2}$  \hfill\vrule}%
\cr
\noalign{\vskip2pt}
\noalign{\hrule}
\noalign{\vskip2pt}
\noalign{\hrule}
(u,~v,~w) & R^{opt} & Vol(W_{uvw}) & Vol(B(R^{opt})) & \Delta^{opt} \cr
\noalign{\hrule}
(8,3,3) & 1.12838 & 0.10721 & 7.75022 & {\mathbf{2.25901}} \cr
\noalign{\hrule}
(7,3,4) & 1.36005 & 0.16297 & 15.19874 & \mathbf{3.33068} \cr
\noalign{\hrule}
(7,3,5) & 1.88213 & 0.23326 & 55.88945 & \mathbf{8.55728} \cr
\noalign{\hrule}
(5,4,3) & 1.28550 & 0.16596 & 12.34723 & \mathbf{3.71986} \cr
\noalign{\hrule}
(4,5,3) & 1.61692 & 0.21299 & 29.64079 & \mathbf{8.69789} \cr
\noalign{\hrule}}}}
\smallbreak
\normalsize
The cases 2.~i.~e,~f,~g,~h with $O=Q,J,E,K$, respectively, will be not relevant for optimal densities.
\subsection{Case 2.ii.b}
In this case $A_3$ is outer, $a_3=JEQ$ is its polar plane, $A_0$ is outer with polar plane $a_0=CLH$. $O=A_2$ is the ball centre (see Fig.~1).

\begin{equation}
\begin{gathered}
r=\min\{A_2b^2,A_2Q,A_2L \};~ \cosh{(A_2b^2)}=\sqrt{1-\frac{1}{a_{22}}},~
\cosh{(A_2Q)}=\frac{1}{\sqrt{a_{33}}}, \\
R=\max\{A_2C, A_2H, A_2J\};~ \cosh{(A_2C)}={\frac{a_{01} a_{02} - a_{12} a_{00}}{\sqrt{a_{00}a_{22}(a_{11}a_{00}-a_{01}^2)}}}, \\
\cosh{(A_2H)}=\frac{\cot{\frac{\pi}{w}}}{\sin \frac{\pi}{v}\sqrt{a_{00}}}, ~
\cosh{(A_2J)}=\frac{\cot{\frac{\pi}{u}}\cot\frac{\pi}{v}}{\sqrt{a_{00}}}.
\end{gathered} \tag{4.12}
\end{equation}

\medbreak
\scriptsize
\centerline{\vbox{
\halign{\strut\vrule~\hfil $#$ \hfil~\vrule
&\quad \hfil $#$ \hfil~\vrule
&\quad \hfil $#$ \hfil\quad\vrule
&\quad \hfil $#$ \hfil\quad\vrule
&\quad \hfil $#$ \hfil\quad\vrule
\cr
\noalign{\hrule}
\noalign{\vskip2pt}
\multispan5{\strut\vrule\hfill\bf Table 2.ii.b. Packing, $O=A_2$, $\frac{1}{u}+\frac{1}{v} < \frac{1}{2},~\frac{1}{v}+\frac{1}{w} < \frac{1}{2}$  \hfill\vrule}%
\cr
\noalign{\vskip2pt}
\noalign{\hrule}
\noalign{\vskip2pt}
\noalign{\hrule}
(u,~v,~w) & r^{opt} & Vol(W_{uvw}) & Vol(B(r^{opt})) & \delta^{opt} \cr
\noalign{\hrule}
(7,3,7) & 0.92836 & 0.38325 & 3.97899 & \mathbf{0.37080} \cr
\noalign{\hrule}
(7,3,8) & 0.94156 & 0.41326 & 4.17152 & \mathbf{0.36051} \cr
\noalign{\hrule}
(7,3, w \to \infty) & \mathrm{arch}\Big(\frac{\sqrt{5-4 \cos^2(\pi/7)}}{\sqrt{1-\cos^2(\pi/7)}} \Big) \approx 0.98513
& 0.49195 & 4.85782 & 0.35267 \cr
\noalign{\hrule}
(5,4,5) & 0.91604 & 0.46190& 3.80539 & \mathbf{0.41193} \cr
\noalign{\hrule}
5,4, w \to \infty) & \mathrm{arch}\Big(\frac{(3-2 \cos^2(\pi/5))}{\sqrt{2(1-\cos^2(\pi/5})}\Big) \approx 1.01789
& 0.59404 & 5.42887 & \mathbf{0.45694} \cr
\noalign{\hrule}}}}
\smallbreak
\scriptsize
\centerline{\vbox{
\halign{\strut\vrule~\hfil $#$ \hfil~\vrule
&\quad \hfil $#$ \hfil~\vrule
&\quad \hfil $#$ \hfil\quad\vrule
&\quad \hfil $#$ \hfil\quad\vrule
&\quad \hfil $#$ \hfil\quad\vrule
\cr
\noalign{\hrule}
\noalign{\vskip2pt}
\multispan5{\strut\vrule\hfill\bf Table 2.ii.b. Covering, $O=A_2$, $\frac{1}{u}+\frac{1}{v} < \frac{1}{2},~\frac{1}{v}+\frac{1}{w} < \frac{1}{2}$  \hfill\vrule}%
\cr
\noalign{\vskip2pt}
\noalign{\hrule}
\noalign{\vskip2pt}
\noalign{\hrule}
(u,~v,~w) & R^{opt} & Vol(W_{uvw}) & Vol(B(R^{opt})) & \Delta^{opt} \cr
\noalign{\hrule}
(7,3,7) & 2.28239 & 0.38325 & 136.50395 & \mathbf{12.72062} \cr
\noalign{\hrule}
(7,3,8) & 2.16470 & 0.41326 & 105.59967 & {9.12600} \cr
\noalign{\hrule}
(7,3,9) & 2.16413 & 0.43171 & 105.46662 & \mathbf{8.72492} \cr
\noalign{\hrule}
(5,4,5) & 1.83500 & 0.46190 & 50.08711 & \mathbf{5.42186} \cr
\noalign{\hrule}
(5,4,6) & 1.79568 & 0.50747 & 45.66776 & \mathbf{4.49957} \cr
\noalign{\hrule}}}}
\smallbreak
\normalsize
The cases 2.~ii.~e,~f,~g,~h with $O=Q,J,E,K,$ respectively are not relevant for optimal densities.
\normalsize
\subsection{Case 2.s.ii.b}
$O=A_2$ is the ball centre (see Fig.~1) and
$u=w$. The density is related to the "half orthoscheme".

\begin{equation}
\begin{gathered}
r=\min\{A_2b^2,A_2Q, A_2F_{12}\};~ \cosh{(A_2b^2)}=\sqrt{1-\frac{1}{a_{22}}},~ 
\cosh{(A_2Q)}=\frac{1}{\sqrt{a_{33}}}, \\
\cosh{(A_2F_{12})}=\sqrt{\frac{1}{2}+\frac{\cos\frac{\pi}{v}}{2\sin^2\frac{\pi}{u}}},~ 
R=\max\{A_2 J,A_2F_{03}\};\\
\cosh{(A_2J)}=\frac{\cot{\frac{\pi}{u}}\cot\frac{\pi}{v}}{\sqrt{a_{00}}},~
\cosh{(A_2F_{03})}=\cot\frac{\pi}{u}\sqrt{\frac{\cos\frac{\pi}{v}+\sin^2\frac{\pi}{u}}
{2(1-\cos\frac{\pi}{v})}}.
\end{gathered} \tag{4.13}
\end{equation}

\scriptsize
\centerline{\vbox{
\halign{\strut\vrule~\hfil $#$ \hfil~\vrule
&\quad \hfil $#$ \hfil~\vrule
&\quad \hfil $#$ \hfil\quad\vrule
&\quad \hfil $#$ \hfil\quad\vrule
&\quad \hfil $#$ \hfil\quad\vrule
\cr
\noalign{\hrule}
\noalign{\vskip2pt}
\multispan5{\strut\vrule\hfill\bf Table 2.s.ii.b. Packing, $O=A_2$,  $\frac{1}{u}+\frac{1}{v}<\frac{1}{2},~\frac{1}{v}+\frac{1}{w}<\frac{1}{2}$  \hfill\vrule}%
\cr
\noalign{\vskip2pt}
\noalign{\hrule}
\noalign{\vskip2pt}
\noalign{\hrule}
(u,~v,~w) & r^{opt} & Vol(W_{uvw}) & Vol(B(r^{opt})) & \delta^{opt} \cr
\noalign{\hrule}
(3,v \to \infty ,3) & \mathrm{arch}(\frac{\sqrt{42}}{6}) \approx 0.39768
& 0.44446 & 0.27191 & \mathbf{0.10196} \cr
\noalign{\hrule}
(4,v \to \infty ,4) & \mathrm{arch}(\frac{\sqrt{6}}{2}) \approx 0.65848
& 0.63434 & 1.30405 & \mathbf{0.25697} \cr
\noalign{\hrule}
(5,5,5) & 0.74746 & 0.57271 & 1.95548 & \mathbf{0.34144} \cr
\noalign{\hrule}
(5,v \to \infty ,3, 5) & \begin{gathered} \mathrm{arch}\Bigg( \sqrt{2-\cos\Big(\frac{\pi}{5}\Big)}  \Bigg) \approx \\ \approx 0.55832 \end{gathered}
& 0.73015 & 0.77585 & 0.10626 \cr
\noalign{\hrule}
(6,4,6) & 0.78340 & 0.55557 & 2.27605 & \mathbf{0.34140} \cr
\noalign{\hrule}
(6,v \to \infty ,6) & \mathrm{arch}(\frac{{\sqrt{5}}}{2}) \approx 0.48121
& 0.78465 & 0.48890 & 0.05192 \cr
\noalign{\hrule}}}}
\smallbreak
\scriptsize
\centerline{\vbox{
\halign{\strut\vrule~\hfil $#$ \hfil~\vrule
&\quad \hfil $#$ \hfil~\vrule
&\quad \hfil $#$ \hfil\quad\vrule
&\quad \hfil $#$ \hfil\quad\vrule
&\quad \hfil $#$ \hfil\quad\vrule
\cr
\noalign{\hrule}
\noalign{\vskip2pt}
\multispan5{\strut\vrule\hfill\bf Table 2.s.ii.b. Covering, $O=A_2$,  $\frac{1}{u}+\frac{1}{v}<\frac{1}{2},~\frac{1}{v}+\frac{1}{w}<\frac{1}{2}$  \hfill\vrule}%
\cr
\noalign{\vskip2pt}
\noalign{\hrule}
\noalign{\vskip2pt}
\noalign{\hrule}
(u,~v,~w) & R^{opt} & Vol(W_{uvw}) & Vol(B(R^{opt})) & \Delta^{opt} \cr
\noalign{\hrule}
(3,8,3) & 1.49263 & 0.32610 & 21.63099 & \mathbf{11.05541} \cr
\noalign{\hrule}
(4,5,4) & 1.43911 & 0.43062 & 18.80212 & \mathbf{5.45785} \cr
\noalign{\hrule}
(5,4,5) & 1.40493 & 0.46190 & 17.16523 & \mathbf{3.71622} \cr
\noalign{\hrule}
(6,4,6) & 1.47083 & 0.55557 & 20.43783 & \mathbf{3.06559} \cr
\noalign{\hrule}
\noalign{\hrule}}}}
\smallbreak
\normalsize
\subsection{Case 2.s.ii.c}
$O=F_{03}$ is the ball centre (see Fig.~1) and
$u=w$. 

\begin{equation}
\begin{gathered}
r=\min\{F_{03}b^0, F_{03}J\};~ \cosh{(F_{03} b^0)}=\sqrt{1+\frac{\cos\frac{\pi}{v}-\sin^2\frac{\pi}{u}}{2(1-\cos\frac{\pi}{v})}}, \\
\cosh{(F_{03}J)}=\frac{\sin\frac{\pi}{v}}{\sqrt{2Ba_{33}(a_{33}+a_{03})}},~
R=\max\{F_{03}A_2, F_{03}Q \};\\
\cosh{(F_{03}A_2)}=\cot\frac{\pi}{u}\sqrt{\frac{\cos\frac{\pi}{v}+\sin^2\frac{\pi}{u}}
{2(1-\cos\frac{\pi}{v})}}, ~ \cosh{(F_{03}Q)}=\frac{\cot\frac{\pi}{u}\cos\frac{\pi}{v}}{\sqrt{2Ba_{33}(a_{33}+a_{03})}}.\\
\end{gathered} \tag{4.14}
\end{equation}

\scriptsize
\centerline{\vbox{
\halign{\strut\vrule~\hfil $#$ \hfil~\vrule
&\quad \hfil $#$ \hfil~\vrule
&\quad \hfil $#$ \hfil\quad\vrule
&\quad \hfil $#$ \hfil\quad\vrule
&\quad \hfil $#$ \hfil\quad\vrule
\cr
\noalign{\hrule}
\noalign{\vskip2pt}
\multispan5{\strut\vrule\hfill\bf Table 2.s.ii.c. Packing, $O=F_{03}$,  $\frac{1}{u}+\frac{1}{v}<\frac{1}{2},~\frac{1}{v}+\frac{1}{w}<\frac{1}{2}$  \hfill\vrule}%
\cr
\noalign{\vskip2pt}
\noalign{\hrule}
\noalign{\vskip2pt}
\noalign{\hrule}
(u,~v,~w) & r^{opt} & Vol(W_{uvw}) & Vol(B(r^{opt})) & \delta^{opt} \cr
\noalign{\hrule}
(3,7,3) & 0.78871 & 0.27899 & 2.32647 & \mathbf{0.59564} \cr
\noalign{\hrule}
(3,8,3) & 0.72041 & 0.32610 & 1.73696 & 0.33290 \cr
\noalign{\hrule}
(4,5,4) & 0.80846 & 0.43062 & 2.52145 & 0.58554 \cr
\noalign{\hrule}
(5,4,5) & 0.72146 & 0.46190 & 1.74508 & 0.47226 \cr
\noalign{\hrule}
(6,4,6) & 0.69217 & 0.55557 & 1.52838 & 0.34388 \cr
\noalign{\hrule}
\noalign{\hrule}}}}
\smallbreak
\scriptsize
\centerline{\vbox{
\halign{\strut\vrule~\hfil $#$ \hfil~\vrule
&\quad \hfil $#$ \hfil~\vrule
&\quad \hfil $#$ \hfil\quad\vrule
&\quad \hfil $#$ \hfil\quad\vrule
&\quad \hfil $#$ \hfil\quad\vrule
\cr
\noalign{\hrule}
\noalign{\vskip2pt}
\multispan5{\strut\vrule\hfill\bf Table 2.s.ii.c. Covering, $O=F_{03}$,  $\frac{1}{u}+\frac{1}{v}<\frac{1}{2},~\frac{1}{v}+\frac{1}{w}<\frac{1}{2}$  \hfill\vrule}%
\cr
\noalign{\vskip2pt}
\noalign{\hrule}
\noalign{\vskip2pt}
\noalign{\hrule}
(u,~v,~w) & R^{opt} & Vol(W_{uvw}) & Vol(B(R^{opt})) & \Delta^{opt} \cr
\noalign{\hrule}
(3,8,3) & 1.26607 & 0.32610 & 11.68147 & \mathbf{2.23886} \cr
\noalign{\hrule}
(4,5,4) & 1.22646 & 0.43062 & 10.41484 & 2.41856 \cr
\noalign{\hrule}
(5,4,5) & 1.28483 & 0.46190 & 12.32350 & 3.33500 \cr
\noalign{\hrule}
(5,5,5) & 1.51882 & 0.57271 & 23.14218 & 4.04082 \cr
\noalign{\hrule}
(6,4,6) & 1.43252 & 0.55557 & 18.47661 & 4.15712 \cr
\noalign{\hrule}}}}

\normalsize
\subsection{Case 2.s.ii.d}
$O=F_{12}$ is the ball centre (see Fig.~1) and
$u=w$. 
\begin{equation}
\begin{gathered}
r=\min\{F_{12}b^1, F_{12}a_3\};~
\cosh{(F_{12}b^1)}=\sqrt{1+\frac{\cos\frac{\pi}{v}-\sin^2\frac{\pi}{u}}{2}}, \\
\cosh{(F_{12} a_3)}=\sqrt{1+\frac{\cos^2\frac{\pi}{u}
(\cos\frac{\pi}{v}+\sin^2\frac{\pi}{u})}{2(\cos^2\frac{\pi}{v}-\sin^2\frac{\pi}{u})}},~ R=\max\{F_{12}J, F_{12}Q \}; \\
\cosh{(F_{12}J)}=\frac{\cos\frac{\pi}{u}(1+\cos\frac{\pi}{v})}{\sin\frac{\pi}{v}\sqrt{2a_{33}(\cos\frac{\pi}{v}+\sin^2\frac{\pi}{u})}},~
\cosh{(F_{12}Q)}=\sqrt{\frac{\cos\frac{\pi}{v}+\sin^2\frac{\pi}{u}}{{2a_{33}\sin^2\frac{\pi}{u}}}}.\\
\end{gathered} \tag{4.15}
\end{equation}
This case is not relevant.
\subsection{Case 2.s.ii.e}
$O=Q$ is the ball centre (see Fig.~1) and
$u=w$. 

\begin{equation}
\begin{gathered}
r=\min\{QA_2,QE,\frac{1}{2}QC\};~
\cosh{(QA_2)}=\frac{1}{\sqrt{a_{33}}},~\cosh{(QE)}=\frac{\cos\frac{\pi}{v}}{\sin\frac{\pi}{u}},\\
\cosh{\Big(\frac{1}{2}QC\Big)}=\sqrt{\frac{1}{2}+\frac{a_{12}}{2a_{22}a_{33}}},~ R=\max\{QF_{03}, QF_{12}\}; \\
\cosh{(QF_{03})}=\frac{\cot\frac{\pi}{u}\cos\frac{\pi}{v}}{\sqrt{2Ba_{33}(a_{33}+a_{03})}}, ~
\cosh{(QF_{12})}=\sqrt{\frac{\cos\frac{\pi}{v}+\sin^2\frac{\pi}{u}}{{2a_{33}\sin^2\frac{\pi}{u}}}}.\\
\end{gathered} \tag{4.16}
\end{equation}
\medbreak
\scriptsize
\centerline{\vbox{
\halign{\strut\vrule~\hfil $#$ \hfil~\vrule
&\quad \hfil $#$ \hfil~\vrule
&\quad \hfil $#$ \hfil\quad\vrule
&\quad \hfil $#$ \hfil\quad\vrule
&\quad \hfil $#$ \hfil\quad\vrule
\cr
\noalign{\hrule}
\noalign{\vskip2pt}
\multispan5{\strut\vrule\hfill\bf Table 2.s.ii.e. Packing, $O=Q$,  $\frac{1}{u}+\frac{1}{v}<\frac{1}{2},~\frac{1}{v}+\frac{1}{w}<\frac{1}{2}$  \hfill\vrule}%
\cr
\noalign{\vskip2pt}
\noalign{\hrule}
\noalign{\vskip2pt}
\noalign{\hrule}
(u,~v,~w) & r^{opt} & Vol(W_{uvw}) & Vol(B(r^{opt})) & \delta^{opt} \cr
\noalign{\hrule}
(3,v \to \infty ,3) & \mathrm{arch}(\frac{2\sqrt{3}}{3}) \approx 0.54931
& 0.44446 & 0.73740 & {0.27652} \cr
\noalign{\hrule}
(4,8,4) & 0.76429 & 0.56369 & 2.10109 & \mathbf{0.46592} \cr
\noalign{\hrule}
(4,9,4) & 0.73883 & 0.57923 & 1.88365 & 0.40650\cr
\noalign{\hrule}
(4,v \to \infty ,4) & \mathrm{arch}(\frac{\sqrt{6}}{2}) \approx 0.65848
& 0.63434 & 1.30405 & 0.25697 \cr
\noalign{\hrule}
(5,5,5) & 0.77537 & 0.57271 & 2.20130 & 0.38436 \cr
\noalign{\hrule}
(5,v \to \infty, 5) & \begin{gathered}
\mathrm{arch}\Bigg( \frac{1}{2}\sqrt{2-\cos\Big(\frac{\pi}{5}\Big)^2}  \Bigg) \approx \\ \approx 0.55832 \end{gathered}
& 0.73015 & 0.77585 & 0.10626 \cr
\noalign{\hrule}
(6,4,6) & 0.78340 & 0.55557 & 2.27605 & 0.34140 \cr
\noalign{\hrule}}}}
\medbreak
\scriptsize
\centerline{\vbox{
\halign{\strut\vrule~\hfil $#$ \hfil~\vrule
&\quad \hfil $#$ \hfil~\vrule
&\quad \hfil $#$ \hfil\quad\vrule
&\quad \hfil $#$ \hfil\quad\vrule
&\quad \hfil $#$ \hfil\quad\vrule
\cr
\noalign{\hrule}
\noalign{\vskip2pt}
\multispan5{\strut\vrule\hfill\bf Table 2.s.ii.e. Covering, $O=Q$,  $\frac{1}{u}+\frac{1}{v}<\frac{1}{2},~\frac{1}{v}+\frac{1}{w}<\frac{1}{2}$  \hfill\vrule}%
\cr
\noalign{\vskip2pt}
\noalign{\hrule}
\noalign{\vskip2pt}
\noalign{\hrule}
(u,~v,~w) & R^{opt} & Vol(W_{uvw}) & Vol(B(R^{opt})) & \Delta^{opt} \cr
\noalign{\hrule}
(3,8,3) & 1.17362 & 0.32610 & 8.90096 & 4.54920 \cr
\noalign{\hrule}
(3,9,3) & 1.23555 & 0.35444 & 10.69512 & 5.02911 \cr
\noalign{\hrule}
(4,5,4) & 1.22646 & 0.43062 & 10.41484 & 3.02321 \cr
\noalign{\hrule}
(4,8,4) & 1.33092 & 0.56369 & 14.02565 & 3.11023 \cr
\noalign{\hrule}
(5,4,5) & 1.28483 & 0.46190 & 12.32350 & 2.66800 \cr
\noalign{\hrule}
(5,5,5) & 1.42880 & 0.57271 & 18.29491 & 3.19444 \cr
\noalign{\hrule}
(6,4,6) & 1.40674 & 0.55557 & 17.24885 & \mathbf{2.58725} \cr
\noalign{\hrule}
(6,5,6) & 1.63628 & 0.64650 & 31.09468 & 4.00806 \cr
\noalign{\hrule}
(6,6,6) & 1.83634 & 0.69130 & 50.24429 & 6.05670 \cr
\noalign{\hrule}}}}
\normalsize
The cases 2.~s.~ii.~f,~g,~h with $O=J,E,K,$ respectively are not relevant for optimal densities.
\medbreak
\bibliographystyle{amsplain}

\end{document}